\documentclass[11pt,a4paper]{amsart}
\usepackage{amssymb, amscd, amsmath, enumerate,verbatim}
\usepackage{pb-diagram, amscd}

\numberwithin{equation}{subsection}
\newtheorem{theorem}{Theorem}[subsection]

\newtheorem{proposition}[theorem]{Proposition}
\newtheorem{corollary}[theorem]{Corollary}

\theoremstyle{definition}

\newtheorem{definition}[theorem]{Definition}
\newtheorem{remark}[theorem]{Remark}
\newtheorem{example}[theorem]{Example}
\newtheorem{examples}[theorem]{Examples}

\newcommand{\mk}{\mathbf{k}}
\newcommand{\id}{\sf id}
\newcommand{\D}{\mathcal{D}}
\newcommand{\A}{\mathcal{A}}

\newcommand{\C}{\mathcal{C}}
\newcommand{\E}{\mathcal{E}}
\newcommand{\X}{\mathcal{X}}

\newcommand{\Sets}{\Delta^o\mathbf{Sets}}
\newcommand{\M}{\mathcal{M}}

\newcommand{\Top}{\mathbf{Top}}

\newcommand{\Ab}{\mathcal{A}}
\newcommand{\Chains}[1]{\mathbf{C}_*(#1)}
\newcommand{\Chainsd}[1]{\mathbf{C}_{**}(#1)}
\newcommand{\Cochainsd}[1]{\mathbf{C}^{**}(#1)}
\newcommand{\Cochains}[1]{\mathbf{C}^*(#1)}
\newcommand{\Chainsab}{\Chains{\Ab}}
\newcommand{\sMod}{\mathbf{\Sigma Mod}}

\newcommand{\Symon}{\mathbf{SyMon}}

\newcommand{\Comon}{\mathbf{CoMon}}

\newcommand{\MOp}{\mathbf{MOp}}
\newcommand{\Op}{\mathbf{Op}}

\newcommand{\Hom}{\mathbf{Cat}}

\begin{document}
\begin{large}

\title[Monoidal functors, acyclic models and chain operads]{Monoidal functors, acyclic models and chain operads}

\author[F. Guill\'{e}n]{F. Guill\'{e}n Santos}
\author[V. Navarro]{V. Navarro}
\address[F. Guill{\'e}n Santos and V. Navarro]{ Departament
d'\`{A}lgebra i Geometria\\  Universitat de Barcelona\\ Gran Via
585, 08007 Barcelona (Spain)}
\author[P.Pascual]{P. Pascual}
\author[A. Roig]{A. Roig}
\address[P. Pascual and A. Roig]{ Departament de Matem\`{a}tica Aplicada
I\\ Universitat Polit\`{e}cnica de Catalunya\\Diagonal 647, 08028
Barcelona (Spain). }
\email{fguillen@ub.edu\\pere.pascual@upc.edu\\agustin.roig@upc.edu
}

\footnotetext[1]{Partially supported by projects DGCYT
BFM2003-60063 and BFM2003-06001}

\date{\today}

\maketitle

\begin{abstract}
We prove that for a topological operad $P$ the operad of oriented
cubical chains, $C^{ord}_\ast(P)$, and the operad of singular
chains, $S_\ast(P)$, are weakly equivalent. As a consequence,
$C^{ord}_\ast(P;\mathbb{Q})$ is formal if and only if
$S_\ast(P;\mathbb{Q})$ is formal, thus linking together some
formality results that are spread out in the literature. The proof
is based on an acyclic models theorem for monoidal functors. We
give different variants of the acyclic models theorem and apply
the contravariant case to study the cohomology theories for
simplicial sets defined by $R$-simplicial differential graded
algebras.
\end{abstract}


\section{Introduction}

Since its introduction by S. Eilenberg and S. MacLane in
\cite{EM1}, the acyclic models theorem has been a powerful
technique in algebraic topology and homological algebra. It says
that, if $K_* , L_* : \A \longrightarrow \Chains{\mathbb Z}$ are
functors from a category with models $\A$ to the category of
non-negative chain complexes of abelian groups $\Chains{\mathbb
Z}$, such that $K_\ast$ is {\em representable} and $L_\ast$ is
{\em acyclic on models}, then any morphism $f_{-1}: H_0(K_{*})
\longrightarrow H_0(L_{*})$ extends to a morphism $f_{*}: K_{*}
\longrightarrow L_{*}$, and this extension is unique up to
homotopy. In particular, if both $K_\ast,L_\ast$, are
representable and acyclic on models, then any isomorphism
$f_{-1}:H_0(K_{*}) \cong H_0(L_{*})$ extends to a homotopy
equivalence $f_{*}: K_{*} \simeq L_{*}$, unique up to homotopy,
thus defining an isomorphism $H{f_*}: H(K_*) \longrightarrow
H(L_*)$. In some elementary presentations, the representability
hypothesis for $K_\ast$ is replaced by asking it to be free on
models, see for instance \cite{D}.

In \cite{EM1}, Eilenberg-MacLane apply their theorem to compare
singular chains $S_\ast(X)$ and normalized cubical chains
$C_\ast(X)$: for $p\geq 0$, the abelian group $S_p(X)$ is freely
generated by continuous maps $\Delta^p\longrightarrow X$, and for
cubical chains one takes first the free group $Q_p(X)$ of
continuous maps $I^p\longrightarrow X$, and defines $C_p(X)$ as
the quotient of $Q_p(X)$ by the degenerated cubical chains
(observe that $C_\ast(X)$ is denoted by $Q^N_\ast(X)$ in
\cite{EM1}). While $S_\ast(X)$ is free on models, $C_\ast(X)$ is
not. Nevertheless, $C_\ast(X)$ is representable, so acyclic models
allows us to extend the identification $S_0(X)=C_0(X)$ to a
natural homotopy equivalence $S_\ast(X)\longrightarrow C_\ast(X)$.
In this way, the two complexes define a homology functor
$H_\ast(X;\mathbb{Z})$ from the category $\Top$ of topological
spaces to the category of graded abelian groups, up to canonical
isomorphism.

Moreover, acyclic models permit to define Eilenberg-Zilber
equivalences in both theories, which together with the K{\"u}nneth
theorem endow the homology functor with morphisms $H_\ast
(X;{\mathbb Z})\otimes H_\ast (Y;{\mathbb Z})\longrightarrow
H_\ast(X\times Y;{\mathbb Z})$, to which we will refer as
K{\"u}nneth morphisms, which are independent of the theory we use,
up to isomorphism. In categorical terms, the K{\"u}nneth morphisms
define a monoidal structure on the functor $H_\ast(\
;\mathbb{Z})$, so we can say that the monoidal functor of homology
is well defined up to isomorphism.

Since the appearance of \cite{EM1}, there have been many
variations and generalizations of the acyclic models theorem,
according to the field of expected applications; see, for
instance, \cite{GMo}, \cite{BB}, \cite{DMO}, \cite {K1} and
\cite{B96}, to cite a few. When dealing with products, acyclic
models allow us to produce a whole family of higher homotopies
from which we can deduce the multiplicativity of morphisms in
cohomology. For example, the natural integration morphism
$\int:S^\ast(X)\longrightarrow \Omega^\ast(X)$ is not
multiplicative, but the higher homotopies defined by acyclic
models define a \emph{dash} morphism which induces a
multiplicative morphism between the singular cohomology and the De
Rham cohomology of $X$, (see \cite{BG}, \cite{Mu}).

For us, the most relevant version of acyclic models is that of M.
Barr in \cite{B96}, also exposed in his book \cite{B02}. The
Eilenberg-MacLane theorem gives morphisms $K_\ast\longrightarrow
L_\ast$ well defined up to homotopy, so that in homology they
define morphisms up to isomorphism. In \cite{B96}, Barr looks for
uniquely defined morphisms in a conveniently localized category of
complexes. In this way, Barr's version of acyclic models
eliminates the indeterminacy up to homotopy fund in the category
of complexes remaining at the complex level and, most remarkably,
it allows to look for extensions defined by a chain of true
morphims $K_\ast\rightarrow\cdot\leftarrow\dots \rightarrow
L_\ast$ in either direction.

Since the functors $S_\ast,C_\ast$ are monoidal functors, they
transform non symmetric operads into non symmetric operads, so we
have the non symmetric {\em dg} operads $S_\ast(P)$ and
$C_\ast(P)$ associated to a given topological operad $P$. Although
one expects these two \emph{dg} operads to be homotopy equivalent,
the known acyclic models theorems does not apply, because the
transformation $S_\ast\longrightarrow C_\ast$ deduced above is not
necessarily a \emph{monoidal} natural transformation of
\emph{monoidal} functors and, consequently, we do not necessarily
obtain an {\em operad} morphism between $S_\ast(P)$ and
$C_\ast(P)$. This situation is reminiscent of what happens with
products in the De Rham theorem cited above.

Our main result, theorem \ref{representable=cofibrante}, is an
acyclic models theorem for monoidal functors from a monoidal
category $\C$ to the monoidal category  $\Chains{\mathbb{Z}}$. We
likewise establish several variations of this result, which cover
the symmetric monoidal and the contravariant monoidal settings. As
a consequence of our results we prove that $S_\ast$ and $C_\ast$
are weakly equivalent as monoidal functors, so it follows that,
for a given topological operad $P$, the non symmetric \emph{dg}
operads $S_\ast(P)$ and $C_\ast(P)$ are weakly equivalent. We
deduce also a multiplicative De Rham comparison result, without
reference to \emph{dash} morphisms.

Theorem \ref{representable=cofibrante} is stated and proved by
adapting to the monoidal setting the categorical framework
proposed by Barr in \cite{B96}. Thus, if $\C$ is a monoidal
category, $K_\ast, L_\ast$ are monoidal functors from $\C$ to
$\Chains{\mathbb{Z}}$, and $H_0(K_\ast)\longrightarrow
H_0(L_\ast)$ is a monoidal morphism, theorem
\ref{representable=cofibrante} establishes sufficient conditions
on $K_\ast , L_\ast$ to extend $f$ to a morphism
$f:K_\ast\longrightarrow L_\ast$ defined up to weak equivalence,
that is, a morphism in the localized category obtained by
inverting a given class of weak equivalences.

The main tools in the statement and proof of theorem
\ref{representable=cofibrante} are the models of $\C$ and the weak
equivalences. The models are given by a suitable monoidal cotriple
$\mathbf G$ defined in $\C$. In the classical situation the
cotriple $\mathbf G$ is a model-induced cotriple. For instance, in
$\Top$ one can take ${\mathbf G}(X) = \sqcup \Delta^n$, where the
disjoint union is taken over all $n\geq 0$ and all continuous maps
$\Delta^n\longrightarrow X$. In the paper, we use a monoidal
version of this cotriple, while when dealing with symmetric
monoidal functors from $\Top$ to $\Chains{\mathbb Z}$ we use the
cotriple introduced by Kleisli in \cite{K2}. As for the class of
morphisms $\Sigma$, in general it is associated to an acyclic
class of complexes (see \cite{B96}), and in the applications it is
a class in between  objectwise homotopy equivalences and
quasi-isomorphisms.

We apply theorem \ref{representable=cofibrante} and its variations
in two different contexts.

First, taking $\C=\Top$, we prove that the singular and the
oriented normalized cubical chain functors are weakly homotopy
equivalent (see \ref{weakformorphisms} for the precise definition)
as symmetric monoidal functors, so their extension to topological
operads give weak-equivalent \emph{dg} operads, (see theorem
\ref{unicidadoperad}). As a consequence, for a topological operad
$P$, the \emph{dg} operads $S_\ast(P)$ and $C_\ast^{ord}(P)$ are
weakly equivalent, so the formality of $P$ does not depend on the
chosen chain functor. This enables us to relate two different
versions of Kontsevich formality theorem for the little cube
operad, the one proved by Kontsevich in \cite{Ko} for oriented
cubical chains and the one by Tamarkin for singular chains (see
\cite{T}), and to give a more precise statement of Deligne's
conjecture.

The second application refers to a uniqueness result for
cohomology theories in the sense of Cartan (see \cite{C}). We
obtain a uniqueness result in the category of \emph{dg} algebras
localized with respect to the class of weak equivalences (see
theorem \ref{unicidadderham} for the precise statement), from
which follows the main theorems in Cartan's paper. In particular,
for any simplicial set $X$ we deduce that the $\mathbb
Q$-differential graded algebras of singular cochains
$S^\ast(X;\mathbb{Q})$ and of Sullivan polynomial forms
$Su^\ast(X)$ are weakly equivalent, recovering a result of
rational homotopy theory, (see \cite{S}, \cite{FHT}).

We now give an overview of the contents of the different sections.
In section $2$, we recall some basic facts about monoidal
functors, and we set some notations used through the paper and
review some relevant examples. We begin section 3 by reviewing the
Eilenberg-MacLane shuffle product formulas and the
Alexander-Whitney formulas, which relate the tensor product in the
category of complexes with the tensor product in the category of
simplicial objects in an additive category $\A$. We use these
general formulas to state the existence of several monoidal and
comonoidal simple functors. Section  4 is devoted to the statement
and proof of theorem \ref{representable=cofibrante}. Following
Barr, (\cite{B96}), we introduce acyclic classes of complexes, the
associated weak equivalences, and the traces of these classes in
the category of monoidal functors. Then, given a monoidal cotriple
$\mathbf G$ in $\C$ we define the relative notions of $\mathbf
G$-presentable monoidal functors from $\C$ to a category of
complexes $\Chains{\A}$ and of $\mathbf G$-weak equivalence for a
monoidal natural transformations between these functors, and prove
the main theorem.

In section 5, we introduce a suitable monoidal triple $\mathbf G$
in the category $\Top$, we prove that singular and normalizad
cubical cochains functors are monoidal $\mathbf G$-presentable and
$\mathbf G$-acyclic, and, as a first application of our general
result, we establish a comparison theorem between this two
monoidal functors.  In section 6, following the same arguments
used in the proof of theorem \ref{representable=cofibrante} and
corollary \ref{unicidadgeneral}, we obtain theorem
\ref{unicidadsim}, which is an acyclic model theorem for symmetric
monoidal functors, from which we prove a comparison result for the
singular and oriented normalized cubical chains as symmetric
monoidal functors, using a triple of Kleisly. As an application we
obtain a similar result for topological operads and topological
modular operads  in section 7. In particular, we obtain the
formality of $C_*^{ord}(\overline{\M},\Bbb Q)$, the $dg$ modular
operad of cubical chains of the modular operad of
Deligne-Knudsen-Mumford moduli spaces of algebraic curves with
marked points. Section 8 is devoted to obtain an acyclic models
theorem for contravariant monoidal functors, and we apply this
result to obtain a comparison result between the singular and
cubical cochains functors on $\Top$. Finally, in section 9 we
apply the previous result to cohomology theories in the sense of
Cartan, including the comparison of Sullivan polynomial forms and
singular cochains.

We should point out that in this paper we do not work out the
contravariant symmetric case. It will be dealt with elsewhere as
an example of a more systematic study of acyclic models in the
context of descent categories (see \cite{GN}).

\section{The monoidal background}

In this section we recall the basic principles of monoidal
categories and functors, (see \cite{EK} or \cite{ML} for details),
and introduce notations that are to be used further in the paper.
These remarks allow us to define in \ref{lossimples} the simple
functors which are associated with simplicial and cosimplicial
functors with values in a category of complexes.

\subsection{Monoidal categories}
\subsubsection{} A {\em monoidal category} consists of a category
$\C$, a functor $\otimes : \C\times \C\longrightarrow \C$, which
we shall call the product functor, and an object $\mathbf{1}$ of
$\C$, in addition to natural isomorphisms of associativity for
$\otimes$ and unit for $\mathbf{1}$, which are subject to
coherence constraints.

\subsubsection{} A monoidal category $(\C,\otimes, \mathbf{1})$ is
{\em symmetric} if, for any objects $X,Y$ of $\C$, it comes
equipped with a natural commutative isomorphism $\tau_{X,Y}
:X\otimes Y\longrightarrow Y\otimes X$ that satisfies the
commutativity constraints.

\subsubsection{} By an {\em additive monoidal category}, also called
an additive tensor category, we understand an additive category
$\A$ which is monoidal in such a way that the product functor
$\otimes : \A\times \A\longrightarrow \A$ is biadditive.

\begin{examples}
In these examples we set out the monoidal structure of several
categories that will appear later on in this paper.

$1$. Let $\A$ be an additive category. We denote by $\Chains{\A}$
the category of uniformly bounded below chain complexes on $\A$,
that is, there is an $n\in \mathbb{Z}$ such that the objects of
$\Chains{\A}$ are differential graded objects $(C_\ast,d)$, with
differential $d$ of degree $-1$, with $C_i=0$ for $i<n$.

If $\A$ is an additive monoidal category, then the category
$\Chains{\A}$ is an additive monoidal category with the product
given by $(C_\ast\otimes D_\ast)_n = \oplus_{p+q=n} C_p\otimes
D_q$, the differential on $C_p\otimes D_q$ being $d=d_C\otimes
{\id} + (-1)^p {\id}\otimes d_D$, and the unit given by the
complex consisting of the unit of $\A$ concentrated in degree
zero.

If $\A$ is a symmetric monoidal category, then $\Chains{\A}$ is
also symmetric, with the natural commutativity isomorphism
$$
\tau_{C_\ast,D_\ast} : C_\ast \otimes D_\ast \longrightarrow
D_\ast \otimes C_\ast
$$
that includes the signs, i.e., it is defined  by
$$
\tau_{C_\ast,D_\ast} = (-1)^{pq} \tau_{C_p,D_q}, \quad
\mbox{on}\quad C_p\otimes D_q.
$$
If $\A$ is the category of $R$-modules for some ring $R$, we will
denote $\Chains{\A}$ by $\Chains{R}$.

$2$. Let $\A$ be an additive monoidal category. By a double
complex of $\A$ we understand a bigraded object of $\A$,
$C_{\ast\ast}$, with two differentials, $d', d''$, which commute;
with the obvious morphisms, they form a category. Thus, the
category of uniformly bounded below double chain complexes
$\Chainsd{\A}$ is isomorphic to $\Chains{\Chains{\A}}$.

$3$.  Let $\X$ be a category and $\D$ be a monoidal category. In
the functor category $\Hom(\X,\D)$ (also denoted by
$\mathbf{Fun}(\X,\C)$, see \cite{B96}) we define the product
$\boxtimes$ componentwise by
$$
(F\boxtimes G) (X) = F(X)\otimes G(X),
$$
and the unit $\mathbf{1}$ as the constant functor defined by the
unit of $\D$. With this structure $\Hom(\X,\D)$ becomes a monoidal
category, which is symmetric if $\D$ is symmetric.

In particular, if $\Delta$ denotes the simplicial category, the
category of simplicial objects of $\D$, which we denote by
$\Delta^o\D$, is a monoidal category. Its objects are denoted by
$X_\bullet$.

Analogously, the category of cosimplicial objects, denoted by
$\Delta\D$, is a monoidal category, which is symmetric if $\D$ is
symmetric. Its objects will be denoted by $X^\bullet$.

$4$. Observe that if $\C$ is a monoidal category, then the dual
category $\C^o$ is also monoidal. In particular, if $º\A$ is an
additive monoidal category, the category of uniformly bounded
below cochain complexes $\Cochains{\A}$, which is isomorphic to
the dual category of $\Chains{\A^o}$, is a monoidal category,
which is symmetric if $\A$ is symmetric. Likewise the category of
uniformly bounded below double chain complexes $\Cochainsd{\A}$ is
a monoidal category.
\end{examples}

\subsection{Monoidal functors}

\subsubsection{} If $(\C,\otimes, \mathbf{1}), (\D,\otimes,
\mathbf{1'})$ are monoidal categories, a {\em monoidal functor\/}
(sometimes called a lax monoidal functor) is a triple
$$
(F , \kappa ,\eta ) : (\mathcal{C}, \otimes , \mathbf{1})
\longrightarrow (\mathcal{D}, \otimes , \mathbf{1}'),
$$
where $F : \mathcal{C} \longrightarrow \mathcal{D}$ is a covariant
functor, together with a natural morphism of $\mathcal{D}$,
$$ \kappa_{X,Y} : FX \otimes FY \longrightarrow F(X \otimes Y),$$
for all objects $X, Y$ of $\mathcal{C}$, and a morphism of
$\mathcal{D}$, $ \eta : \mathbf{1'} \longrightarrow F\mathbf{1}$,
compatible with the constraints of associativity and unit. We will
refer to $\kappa$ as the {\em K\"{u}nneth morphism} of $F$.

A monoidal functor $F$ is said to be a {\em strong monoidal}
functor if the Kunneth morphisms are isomorphisms.

\subsubsection{} If $\mathcal C$ and $\mathcal D$ are symmetric
monoidal categories, a monoidal functor $F : \mathcal{C}
\longrightarrow \mathcal{D}$ is said to be {\em symmetric} if
$\kappa$ is compatible with the commutativity isomorphisms $\tau$.

\begin{examples}\label{exempleToT}
Let $\A$ be an additive (symmetric) monoidal category.

$1$. Given a double complex $C_{\ast\ast}$ of $\A$ with finite
anti-diagonals, the associated total complex is the complex
$\mbox{Tot}(C_{\ast\ast})$ which in degree $n$ is given by
$\mbox{Tot}_n(C_{\ast\ast}) = \oplus_{p+q=n} C_{pq}$ with
differential defined by $d=d'+(-1)^pd''$ on $C_{pq}$. For bounded
below double complexes the total complex defines a functor
$\mbox{\rm Tot}: \Chainsd{\A}\longrightarrow \Chains{\A}$. It is a
strong (symmetric) monoidal functor.

Similarly, the total complex functor defines a strong (symmetric)
monoidal functor for cochain complexes $\mbox{\rm Tot}:
\Cochainsd{\A}\longrightarrow \Cochains{\A}$.

$2$. The {\em homology functor\/} $H_*: \Chainsab \longrightarrow
\Chainsab$ is a (symmetric) monoidal functor, taking as $\kappa$
the usual K{\"u}nneth morphism
$$\kappa_{X,Y}: H_*(X)\otimes H_*(Y) \longrightarrow H_*(X \otimes Y).$$
\end{examples}

\subsection{} Monoidal functors are stable under composition, so the
category of monoidal categories and monoidal functors,
$\mathbf{Mon}$, is a subcategory of $\Hom$.

\subsubsection{} Let $F,G:\mathcal{C}\longrightarrow \mathcal{D}$
be two monoidal functors. A natural transformation $ \phi : F
\Rightarrow G $ is said to be {\em monoidal\/} if it is compatible
with $\kappa$ and $\eta$.

Monoidal functors between two monoidal categories $\C, \D$,
together with monoidal natural transformations, define a
subcategory of the functor category $\Hom(\C,\D)$ which will be
denoted by ${\mathbf{Mon}(\C,\D)}$.

Moreover, monoidal functors and monoidal natural transformations
are compatible, so they equip the class of monoidal categories,
$\mathbf{Mon}$, with a $2$-category structure.

\subsubsection{} Similarly, one can extend the notions above to the
symmetric setting. Hence, the symmetric monoidal functors between
two symmetric monoidal categories $\C, \D$, together with
symmetric monoidal natural transformations define a subcategory of
the functor category $\Hom(\C,\D)$ which is denoted by
${\Symon(\C,\D)}$.

\subsection{Comonoidal functors} After the definition of monoidal functor
we obtain, by duality, the notion of comonoidal functor between
monoidal categories (called op-monoidal functor in \cite{KS}).

\subsubsection{}
By definition, a {\em comonoidal functor} between $\C$ and $\D$ is
a covariant monoidal functor between the dual categories
$\C^o,\D^o$.

That is to say, a comonoidal functor $F:\C\longrightarrow \D$ is a
covariant functor together with natural morphisms
$$
\kappa^o : F(X\otimes Y) \longrightarrow F(X)\otimes F(Y),
$$
for objects $X,Y$ of $\C$, and a morphism $\eta^o :
F\mathbf{1}\longrightarrow \mathbf{1}'$ of $\D$, satisfying
constraints of associativity and unit.
If $\C,\D$, are symmetric monoidal categories, a comonoidal
functor $F$ is said to be symmetric if $\kappa^o$ is compatible
with the commutativity constraint.


\subsubsection{} The definition of comonoidal natural
transformations is clear. We denote by $\Comon (\C,\D)$ the
category of comonoidal functors and comonoidal natural
transformations. The category $\Comon(\C,\D)$ is the dual category
of $\mathbf{Mon}(\C^o,\D^o)$, so we will focus our attention on
monoidal functors.

\subsubsection{}\label{comon=mon}
We have defined the notions of monoidal and comonoidal functors
for covariant functors between two monoidal categories $\C,\D$. It
is straightforward to define the corresponding notions for
contravariant functors from $\C$ to $\D$. Nevertheless, we should
point out that a contravariant monoidal functor $\C\longrightarrow
\D$ (resp. a contravariant comonoidal functor) is equivalent to a
comonoidal functor (resp. monoidal functor) defined in the dual
category, $\C^o\longrightarrow \D$, so nothing new is gained.

\section{Shuffle and Alexander-Whitney maps}

\subsection{Simplicial chain complexes}\label{simple} Let $\A$ be
an additive monoidal category. We use $\sim : \Delta^o
\A\longrightarrow \Chains{\A}$ to denote the functor which takes a
simplicial object $C_\bullet$ of $\A$ to the chain complex
$\widetilde{C_\bullet}$ given by $(\widetilde{C_\bullet})_n= C_n$
with differential
$d=\partial_0-\partial_1+\dots+(-1)^n\partial_n$.

\subsubsection{Shuffle map}\label{shufflesformulas}
The functor $\sim$ with the shuffle map is monoidal (see
\cite{EM2}): recall that if $C_{\bullet}, D_{\bullet}$ are objects
of $\Delta^o{\A}$, the shuffle product $sh : C_{\ast}\otimes
D_{\ast}\longrightarrow \widetilde{C_{\bullet}\boxtimes
D_{\bullet}}$ is defined in degree $n=p+q$, $C_{p}\otimes
D_{q}\longrightarrow C_{p+q}\otimes D_{p+q}$, by the formula
$$
sh_{pq} = \sum_{(\mu, \nu)} \varepsilon(\mu,\nu) (s_{\nu_q}\cdots
s_{\nu_1}) \otimes (s_{\mu_p}\cdots s_{\mu_1}),
$$
where the sum is taken over all $(p,q)$-shuffles $(\mu, \nu)$ and
$\varepsilon(\mu,\nu)$ is the signature of the associated
permutation.

\subsubsection{Alexander-Whitney map}\label{AWformulas}
The functor $\sim $ with the Alexan\-der-Whitney map is also
comonoidal, (loc. cit.). Recall that if $C_\bullet, D_\bullet$ are
simplicial objects of $\A$, the Alexander-Whitney morphism $AW :
\widetilde{C_\bullet\boxtimes D_\bullet}\longrightarrow  C_\ast
\otimes D_\ast$ is given by morphisms $AW : C_n\otimes
D_n\longrightarrow \oplus_{i=0}^n C_i \otimes D_{n-i}$ which are
defined by
$$
AW = \sum_{i=0}^n \widetilde{\partial}^{n-i} \otimes
\partial_0^i,
$$
where $\widetilde{\partial}^{p-i}$ is the last face operator
$\partial_{i+1}\dots\partial_p$.

\subsubsection{Cosimplicial cochain complexes}
Let $\A$ be an additive monoidal category. The category of
cosimplicial objects of $\A$, $\Delta\A$, is a monoidal category.
Let $\sim:\Delta {\A}\longrightarrow \Cochains {\A}$ represent the
functor which takes a cosimplicial object $K^\bullet$ to the
cochain complex $K^\ast$ with differential $d=\partial0+ \dots +
(-1)^n\partial^n$. By dualizing \ref{shufflesformulas} and
\ref{AWformulas}, one obtains comonoidal and monoidal structures,
respectively, for this functor $\sim$. We use the monoidal
structure issuing from the Alexander-Whitney formulas to study
contravariant functors in \S 8.

\subsection{The simple functors}\label{lossimples}
When $\A$ is an additive monoidal category we can use the shuffle
and Alexander-Whitney structures above to put monoidal and
comonoidal structures in the classic simple complex associated to
a simplicial complex of $\A$.

\subsubsection{}\label{eilenbergmac} As $\Chains{\A}$
is an additive monoidal category, the functor $\sim$ defines a
functor
$$
\Delta^o\Chains{\A}\longrightarrow 
\Chainsd {\A}, \eqno{(1)}
$$
which, with the shuffle product, is monoidal. The composition of
this functor with the total functor (which, as pointed out in
example \ref{exempleToT}, is monoidal) is a monoidal functor
$$
s_{EM} : \Delta^o\Chains{\A}\longrightarrow \Chains{\A},
$$
which will be called the {\em Eilenberg-MacLane simple functor}
(or the simple functor, for short).

\subsubsection{} \label{simplechuflas} If $\C, \D$ are monoidal
categories, there is an equivalence of categories $\Delta^o
\mathbf{Mon} (\C,\D) \cong \mathbf{Mon} (\C,\Delta^o \D)$, so if
$\A$ is an additive monoidal category, the composition of this
isomorphism with the Eilenberg-MacLane simple functor is a
functor, also denoted by $s_{EM}$,
$$
s_{EM}: \Delta^o \mathbf{Mon} (\C,\Chains{\A}){\longrightarrow}
\mathbf{Mon} (\C,\Chains{\A}),
$$
which associates a monoidal functor to a simplicial monoidal
functor between $\C$ and $\Chains{\A}$. We will also refer to this
composition as the {\em Eilenberg-MacLane simple functor}.

\subsubsection{} If in \ref{eilenbergmac} we use the comonoidal structure
on $\sim $ that comes from the Alexander-Whitney map, we obtain a
comonoidal functor
$$
s_{AW} : \Delta^o\Chains{\A}\longrightarrow \Chains{\A},
$$
which we call the {\em Alexander-Whitney simple functor}. Thus, in
a way completely analogous to \ref{simplechuflas}, we obtain a
simple functor for comonoidal functors
$$
s_{AW}:\Delta^o \mathbf{CoMon} (\C,\Chains{\A}){\longrightarrow}
\mathbf{CoMon} (\C,\Chains{\A}).
$$

\subsubsection{}\label{AWcontra} By duality, there are similarly
defined simple functors for the categories of monoidal and
comonoidal functors with values in the cochain category
$\Cochains{\A}$. In \S 8 we will use, in the contravariant
setting, the Alexander-Whitney simple functor
$$
s_{AW}:\Delta \mathbf{Mon} (\C^o,\Cochains{\A}){\longrightarrow}
\mathbf{Mon} (\C^o,\Cochains{\A}).
$$
Although $s_{AW}$ denotes the two functors introduced above, it
will be clear from the context which we are referring to.

\subsection{Symmetric simple functor} \label{symsimple} Let $\A$
be an additive symmetric monoidal category. As the shuffle product
is symmetric, (\cite{EM2}), the Eilenberg-MacLane simple functor
$s_{EM}:\Delta^o\Chains{\A}\longrightarrow \Chains{\A}$ is a
symmetric monoidal functor. Thus, if $\C$ is a symmetric monoidal
category, we can follow the reasoning in \ref{simplechuflas} to
deduce the existence of an {\em Eilenberg-MacLane symmetric simple
functor}
$$
s_{EM}:\Delta^o\Symon(\C,\Chains{\A})\longrightarrow
\Symon(\C,\Chains{\A}).
$$

\section{An acyclic models theorem for monoidal functors}
In this section we prove the main technical tool, the acyclic
models theorem for monoidal functors. Our presentation is a
variation of the scheme devised by Barr in \cite{B96} for these
kinds of results. According to Barr, three main ingredients are
required to state an acyclic models theorem: the total complex
associated to a double complex, the acyclic classes of complexes
that define the associated classes of weak equivalences, and the
cotriples constructed from the models. The total complex functor
and its monoidal counterparts were the object of the previous
section, we begin now by introducing acyclic classes.

\subsection{Acyclic classes and weak equivalences}
Let $\A$ be an abelian category. If $C_\ast$ is a chain complex on
$\A$, we denote by $C_\ast[-1]$ the chain complex given by
$(C_\ast [-1])_n=C_{n-1}$ with differential defined by $-d_{n-1}$.

\subsubsection{} Recall the following definition from \cite{B96},
(see also \cite{B02}),

\begin{definition}
A class $\Gamma$ of complexes in $\Chains{\A}$ is called an {\em
acyclic class} if the following conditions are satisfied:
\begin{itemize}
\item[{\sf (AC1)}] The complex $0$ is in $\Gamma$.\item[{\sf
(AC2)}] Stability: $C_\ast$ is in $\Gamma$ if and only if
$C_\ast[-1]$ is in $\Gamma$. \item[{\sf (AC3)}] Let $C_\ast,
D_\ast$ be two chain homotopic complexes. Then $C_\ast$ is in
$\Gamma$ if and only if $D_\ast$ is in $\Gamma$. \item[{\sf
(AC4)}] Every complex in $\Gamma$ is acyclic. \item[{\sf (AC5)}]
If $C_{\ast\ast}$ is a double complex in $\Chainsd{\A}$ all of
whose rows are in $\Gamma$, then the total complex $\rm{Tot}\,
C_{\ast\ast}$ is in $\Gamma$.
\end{itemize}
\end{definition}

By {\sf (AC1)} and {\sf (AC3)} every contractible complex is in
$\Gamma$ and by {\sf (AC4)} all complexes in $\Gamma$ are acyclic,
thus $\Gamma$ is a class between the class of all contractible
complexes and the class of all acyclic complexes. These two
extreme cases are examples of acyclic classes, (see \cite{B96}, \S
4).

\subsubsection{} Given an acyclic class $\Gamma$ in $\Chains{\A}$, we
denote by $\Sigma$ the class of morphisms between chain complexes
of $\A$ whose mapping cone is in $\Gamma$. Morphisms in $\Sigma$
are referred to as {\em weak equivalences} (with respect to
$\Gamma$).

By {\sf (AC1)} and {\sf (AC3)} every homotopy equivalence is in
$\Sigma$ and by $\sf (AC4)$ the mapping cone of a morphism $f$ of
$\Sigma$ is acyclic, thus it follows from the exact sequence of
the mapping cone that any such $f$ is a quasi-isomorphism.

\begin{remark}
Let ${\mathbf K}_\ast (\A)$ be the category of chain complexes up
to homotopy, that is, its objects are complexes of $\A$ and its
morphisms are homotopy classes of morphisms of complexes. It is a
triangulated category. An acyclic class $\Gamma$ determines a
triangulated subcategory of ${\mathbf K}_\ast (\A)$ and, as a
consequence, the class of morphisms $\Sigma$ associated with
$\Gamma$ inherits some properties from the general setup of
triangulated categories, (i.e. it admits a calculus of fractions).
For the sake of simplicity we will follow Barr's treatment and
refer to \cite{B95} for the properties of $\Sigma$ that will be
used.
\end{remark}

\subsubsection{}\label{cadenaweak} We use $\Chains{\A}[\Sigma^{-1}]$ to denote the
localized category of $\Chains{\A}$ with respect to $\Sigma$,
which exists in a suitable universe, (see \cite{GZ}), and is
uniquely determined up to isomorphism.

The class of morphisms $\Sigma$ is stable under composition and
satisfies the 2 out of 3 property (see the proof of proposition
$3.3$ in \cite{B96}), that is, for every pair of morphisms $f,g$
of $\Chains{\A}$ so that $gf$ exists, if {\em two} of $f$, $g$ and
$gf$ are in $\Sigma$, then so is the third. Moreover, $\Sigma $
has a homotopy calculus of fractions (see loc. cit.). In many
cases, such as when $\Sigma$ is the class of quasi-isomorphisms or
the class of homotopy equivalences, $\Sigma$ is a saturated class
of morphisms, in other words, it is precisely the class of
morphisms of $\Chains{\A}$ which become isomorphisms in
$\Chains{\A}[\Sigma^{-1}]$.

Two objects $C_\ast, D_\ast$ of $\Chains{\A}$ are said to be {\em
weakly equivalent} (with respect to $\Gamma$) if there exists a
sequence of morphisms of $ \Chains{\A}$,
$$
C_\ast\longleftarrow C_\ast1\longrightarrow \dots \longleftarrow
C_\ast^m\longrightarrow D_\ast,
$$
which are weak equivalences. Hence, weakly equivalent objects are
isomorphic in $\Chains{\A}[\Sigma^{-1}]$.

\subsubsection{}\label{weakformorphismsone} Let $\C$ be a category and $\A$ an abelian
category. The functor category $\Hom(\C,\A)$ is an abelian
category. We fix a class $\Sigma$ of weak equivalences in
$\Hom(\C,\Chains{\A})$ associated to an acyclic class $\Gamma $ in
$\Chains{\Hom(\C,\A)}=\Hom(\C,\Chains{\A})$.

Let us now assume that $\C$ is a monoidal category and $\A$ is an
abelian monoidal category. As $\mathbf{Mon}(\C, \Chains{\A})$ is a
subcategory of $\Hom(\C,\Chains{\A)}$, $\Sigma$ determines a class
of morphisms in $\mathbf{Mon}(\C, \Chains{\A})$, which will be
represented by the same symbol. Thus, a morphism of
$\mathbf{Mon}(\C,\Chains{\A})$ is a weak equivalence if it is in
$\Sigma$ as a morphism of $\Hom(\C, \Chains{\A})$.

However, note that for two monoidal functors from $\C$ to
$\Chains{\A}$, the weak equivalence relation in $\mathbf{Mon}(\C,
\Chains{\A})$ is not the same as it is in $\Hom(\C, \Chains{\A})$,
since in the first case the intermediate functors in the chain
\ref{cadenaweak} have to be monoidal.

We use $\mathbf{Mon}(\C, \Chains{\A})[\Sigma^{-1}]$ to denote the
category obtained by inverting the weak equivalences in $\Sigma$,
so the natural functor $\mathbf{Mon}(\C,
\Chains{\A})\longrightarrow \mathbf{Mon}(\C,
\Chains{\A})[\Sigma^{-1}]$ transforms weak equivalences to
isomorphisms.

\subsubsection{}\label{arasi}
We note that the class of weak equivalences in
$\mathbf{Mon}(\C,\Chains{\A})$ contains the homotopy equivalences
and is compatible with the functor $s_{EM}$, so it defines a
functor
$$
s_{EM} : (\Delta^o\mathbf{Mon}(\C,
\Chains{\A}))[\widetilde{\Sigma}^{-1}]\longrightarrow
\mathbf{Mon}(\C, \Chains{\A})[\Sigma^{-1}],
$$
where $\widetilde\Sigma$ is the class of morphisms in
$\Delta^o\mathbf{Mon}(\C, \Chains{\A})$ which in each simplicial
degree are in $\Sigma$, as follows from {\sf (AC5)}.

This compatibility is one of the basic properties of acyclic
classes that is needed to prove the acyclic models theorem below.
Instead of Barr's acyclic classes we could work in other settings
where there is such a compatibility (see also \ref{descenso}).

\newpage

\subsection{$\mathbf G$-presentable objects and $\mathbf G$-weak equivalences}
Let $\C$ be a monoidal category, $\A$ an abelian monoidal
category, $\Sigma$ a class of weak equivalences in
$\Hom(\C,\Chains{\A})$ which contains the homotopy equivalences
and is compatible with $s_{EM}$.

\subsubsection{The standard construction}\label{bar}
Recall that if $\X$ is a category, a cotriple (also called a {\em
comonad}, see \cite{ML}) $\mathbf{G}=(G,\varepsilon, \delta)$ in
$\X$ is given by a functor $G:\X \longrightarrow \X$ and natural
transformations $\varepsilon : G\Rightarrow \id$ and $\delta:
G\Rightarrow G2$, which satisfies
$$
\delta G \cdot \delta = G\delta\cdot\delta : G\Rightarrow G3,\quad
\varepsilon G\cdot\delta = \mathbf{1}_G = G\varepsilon\cdot\delta
: G\Rightarrow G,
$$
where the dot denotes the composition of natural transformations,
(see \cite{ML}). Given a cotriple $\mathbf{G}$ in $\X$, every
object $X$ of $\X$ has a functorial augmented simplicial object
associated with it, $B_\bullet (X)$, which will be called the {\em
standard construction} of $\mathbf G$ applied to $X$, (\cite{ML},
chapter VII). In degree $n$, $n\geq 0$, the simplicial object
$B_\bullet(X)$ is given by $G^{n+1}(X)$ with face and degeneracy
transformations given by
$$
\begin{array}{lll}
\partial_i =&  G^i \varepsilon G^{n-i} : G^{n+1}(X)
\longrightarrow G^n(X),&\ 0\leq i\leq n,\\
s_i =& G^i \delta G^{n-i} : G^{n+1}(X) \longrightarrow G^{n+2}(X),
&\ 0\leq i\leq n.
\end{array}
$$
In this way we obtain a simplicial object with an augmentation
defined by $\varepsilon: B_0(X)=G(X) \longrightarrow X$.

By the functoriality of the standard construction we obtain an
augmented simplicial functor $B_\bullet$ in $\Hom(\X,\X)$.

\begin{proposition}\label{GBG}
\begin{itemize}
\item[(1)] The augmented simplicial functor $G(\varepsilon):
G\circ B_\bullet \Rightarrow G $ is contractible. \item[(2)] The
augmented simplicial functor $G(\varepsilon):B_\bullet \circ
G\Rightarrow G$ is contractible.
\end{itemize}
\end{proposition}
\begin{proof}
(1) The face and degeneracy morphisms of $G\circ B_\bullet $ are
given by
$$
\begin{array}{llll}
G(\partial_i^n) =&  G^{i+1} \varepsilon G^{n-i} &=
\partial_{i+1}^{n+1},
&\ 0\leq i\leq n,\\
G(s_i^n) =& G^{i+1}\delta G^{n-i} &= s_{i+1}^{n+1}, &\ 0\leq i\leq
n,
\end{array}
$$
where the latter are the faces and degeneracies of $B_\bullet $.
Now, the extra degeneracy from $B_\bullet $, $s=s_0^{n+1}$, gives
a contraction for $G\circ B_\bullet$.

For (2) use the last extra degeneracy $s^{n+1}_{n+1}$ in $B_n\circ
G$.
\end{proof}

\subsubsection{Monoidal structure in the standard construction}
Let $\mathbf{G}$ be a cotriple in $\mathbf{Mon} (\C, \Chains
{\A})$. We say that $\mathbf G$ is {\em compatible} with $\Sigma$
if $G(\Sigma)\subseteq \Sigma$. We will associate to $\mathbf G$ a
monoidal functor $ B_\ast : \mathbf{Mon}(\C, \Chains{\A})
\longrightarrow \mathbf{Mon}(\C, \Chains{\A})$.

The standard construction defines a functor
$$
B_\bullet : \mathbf{Mon}(\C, \Chains{\A}) \longrightarrow
\Delta^o\mathbf{Mon}(\C, \Chains{\A}).
$$
By composing $B_\bullet$ with the Eilenberg-MacLane simple functor
defined in \ref{simplechuflas}, $s_{EM}$, we obtain a functor
$$
B_\ast = s_{EM}B_\bullet :\mathbf{Mon}(\C, \Chains{\A})
\longrightarrow \mathbf{Mon}(\C, \Chains{\A}).
$$
The natural transformation $\varepsilon$ gives a natural
transformation $\varepsilon : B_\ast\Rightarrow \id$.

\subsubsection{$\mathbf G$-presentable objects} Barr introduces in
\cite{B96} the $\varepsilon$-presentable objects. We will refer to
these objects as $\mathbf G$-presentable objects. Recall its
definition:

\begin{definition}
We say that $K_\ast$ of $\mathbf{Mon}(\C, \Chains{\A})$ is
$\mathbf G$-{\em presentable} (with respect to $\Gamma$) if the
augmentation morphism $\varepsilon : B_\ast(K_\ast)\longrightarrow
K_\ast $ is a weak equivalence.
\end{definition}

As remarked in \ref{arasi}, $\Sigma$ contains the homotopy
equivalences, so if $K_\ast$ is an object of $\mathbf{Mon}(\C,
\Chains{\A})$ such that $\varepsilon_{K_n}$ splits for all $n$,
then $\varepsilon_{K_n}\in \Sigma$, and since $\Sigma$ is
compatible with $s_{EM}$, it follows that $\varepsilon_{K}\in
\Sigma$. Hence we obtain the following result,

\begin{proposition}\label{split}
Let $K_\ast$ be an object of $\mathbf{Mon}(\C, \Chains{\A})$ such
that $\varepsilon_{K_n}$ splits for all $n\in \mathbb{Z}$, that is
to say, for each $n$ there is a natural transformation $\theta_n :
K_n\longrightarrow K_n G$ such that
$\varepsilon_{K_n}\theta_n=\id$. Then, $K_\ast$ is $\mathbf
G$-presentable.
\end{proposition}

\begin{example}\label{GKpresentable}
If $K_\ast$ is an object of $\mathbf{Mon}(\C, \Chains{\A})$, then
$G(K_\ast)$ is $\mathbf G$-presentable, since for each $n$ we can
split $G(K_n)$ by $\delta$.
\end{example}

\subsubsection{G-Weak equivalences} Barr introduces in \cite{B96}
the $\mathbf G$-acyclic objects. More generally we can speak of
objects that are $\mathbf G$-equivalent in the following sense,

\begin{definition}
Let $f: K_\ast\longrightarrow L_\ast$ be a morphism of
$\mathbf{Mon}(\C, \Chains {\A})$. We say that $f$ is a
$B\mathbf{G}$-{\em weak equivalence}, (respectively, a
$\mathbf{G}$-{\em weak equivalence}), with respect to $\Gamma$, if
$B_\ast(f)\in \Sigma$, (respectively, if $G(f)\in \Sigma$).
\end{definition}

\begin{proposition}\label{weakG=weakBG}
If $G(\Sigma)\subseteq \Sigma$, then  a morphism $f:
K_\ast\longrightarrow L_\ast$ of $\mathbf{Mon}(\C, \Chains {\A})$
is a $\mathbf{G}$-weak equivalence if, and only if, it is a
$B\mathbf G$-weak equivalence.
\end{proposition}
\begin{proof}
If $f$ is a $\mathbf G$-weak equivalence then, for all $n\geq 0$,
$G^{n+1}(f)$ is a weak equivalence, by hypothesis. Therefore, by
{\sf (AC5)}, $f$ is a $B\mathbf G$-weak equivalence.

Reciprocally, if we assume that $f$ is a $B\mathbf G$-weak
equivalence, then $G(B(f)): G(B_\ast(K_\ast))\longrightarrow
G(B_\ast(L_\ast))$ is a weak equivalence, by hypothesis. Consider
the commutative diagram
$$
\begin{diagram}
\node{G(B_\ast(K_\ast))}\arrow{s,l}{G(\varepsilon)}\arrow{e,t}{G(B(f))}
\node{G(B_\ast(L_\ast))} \arrow{s,r}{G(\varepsilon)}\\
\node{G(K_\ast)}\arrow{e,b}{G(f)}\node{G(L_\ast)}
\end{diagram}
$$
By proposition \ref{GBG} (2), the two vertical morphisms are
isomorphisms, so the result follows from the 3 out of 2 property
of $\Sigma$.
\end{proof}

Now we recover Barr's definition:

\begin{definition}
An object $K_\ast$ of $\mathbf{Mon}(\C, \Chains{\A})$, in
non-negative degrees, (that is, with $K_p=0$ if $p<0$), is said to
be {\em $\mathbf{G}$-acyclic} (with respect to $\Gamma$) if the
augmentation $K_\ast\longrightarrow H_0(K_\ast)$ is a
$\mathbf{G}$-weak equivalence.
\end{definition}

\subsection{The acyclic models theorem for monoidal functors}
The following result and its corollaries are a variation of the
classical acyclic models theorem in the context of monoidal
functors (see \cite{EM1} and also \cite{BB}, \cite{B96}).

Let $\C$ be a monoidal category, $\A$ an abelian monoidal
category, and $\Sigma$ a class of weak equivalences in
$\mathbf{Mon}(\C, \Chains {\A})$ which contains the homotopy
equivalences and is compatible with $s_{EM}$. Let $\mathbf G$ be a
cotriple in $\mathbf{Mon}(\C, \Chains {\A})$.

\begin{theorem}\label{representable=cofibrante}
Let $K_\ast$ be a $\mathbf G$-presentable object of
$\mathbf{Mon}(\C, \Chains {\A})$ and $s:L_\ast\longrightarrow
M_\ast$ a $\mathbf G$-weak equivalence. Suppose that $\mathbf G$
is compatible with $\Sigma$, (that is, $G(\Sigma)\subseteq
\Sigma$). Then, for any $\alpha : K_\ast\longrightarrow M_\ast$
there exists a unique morphism $\widetilde \alpha :
K_\ast\longrightarrow L_\ast $ in $\mathbf{Mon}(\C,
\Chains{\A})[\Sigma^{-1}]$ such that $\alpha =
s\widetilde{\alpha}$ in this localized category.
\end{theorem}
\begin{proof}
By the naturality of $\varepsilon$, we have a commutative diagram
$$
\begin{diagram}
\node{B_\ast(K_\ast)}\arrow{s,r}{\varepsilon_K}\arrow{e,t}{B(\alpha)}
\node{B_\ast(M_\ast)}\arrow{s,r}{\varepsilon_M}
\node{B_\ast(L_\ast)}\arrow{s,r}{\varepsilon_L}\arrow{w,t}{B(s)}\\
\node{K_\ast}\arrow{e,t}{\alpha} \node{M_\ast}
\node{L_\ast}\arrow{w,t}{s}
\end{diagram}
$$
As $K_\ast$ is $\mathbf G$-presentable, $\varepsilon_K$ is a weak
equivalence.
Moreover, by proposition \ref{weakG=weakBG}, $B(s)$ is also a weak
equivalence. So they are isomorphisms in the localized category
$\mathbf{Mon}(\C, \Chains{\A})[\Sigma^{-1}]$.

In this category we define the morphism
$$
\widetilde{\alpha} = \varepsilon_L B(s)^{-1}B(\alpha
)(\varepsilon_K)^{-1}.
$$
We have $s\widetilde{\alpha} = \alpha$, since by the commutativity
of the diagram above we have $\alpha = \varepsilon_M
B(\alpha)(\varepsilon_K)^{-1}$ and $\varepsilon_M = s\varepsilon_L
B(s)^{-1}$, so it follows that
$$
\alpha = \varepsilon_M B(\alpha)(\varepsilon_K)^{-1} =
s\varepsilon_L B(s)^{-1} B(\alpha)(\varepsilon_K)^{-1} =
s\widetilde{\alpha}.
$$
With respect to uniqueness, assume that $\gamma :K_\ast
\longrightarrow L_\ast$ is another lifting of $\alpha $ in
$\mathbf{Mon}(\C, \Chains{\A})[\Sigma^{-1}]$, so that $\alpha =
s\gamma $. Since the standard construction is functorial and
compatible with weak equivalences because $G(\Sigma)\subseteq
\Sigma$, we have $B(s)B(\widetilde\alpha)=B(\alpha) = B(s)B(\gamma
)$, but $B(s)\in \Sigma$ so $B(\gamma)=B(\widetilde\alpha)$.
Moreover, as $\varepsilon $ is a natural transformation, we have
$\varepsilon_L B(\gamma ) = \gamma \varepsilon_K$ and
$\varepsilon_L B(\widetilde\alpha ) = \widetilde\alpha
\varepsilon_K$, and since $\varepsilon_K\in \Sigma$ we deduce that
$\gamma = \widetilde\alpha$. This ends the proof of the theorem.
\end{proof}

\begin{corollary}
Let $K_\ast,L_\ast$ be objects of $\mathbf{Mon}(\C, \Chains{\A})$
in non-negative degrees. Suppose that $K_\ast$ is $\mathbf
G$-presentable and $L_\ast$ is $\mathbf{G}$-acyclic. Then, any
monoidal natural transformation $H_0(K_\ast)\longrightarrow
H_0(L_\ast)$ has a unique extension to a morphism $K_\ast
\longrightarrow L_\ast$ in $\mathbf{Mon}(\C,
\Chains{\A})[\Sigma^{-1}]$.
\end{corollary}
\begin{proof} It follows from the previous theorem applied to the
diagram
$$
\begin{diagram}
\node[3]{L_\ast}\arrow{s}\\
\node{K_\ast}\arrow{ene,..}\arrow{e}\node{H_0(K_\ast)}
\arrow{e}\node{H_0(L_\ast)}
\end{diagram}
$$
\end{proof}

\begin{corollary}\label{unicidadgeneral}
Let $K_\ast,L_\ast$ be objects of $\mathbf{Mon}(\C, \Chains{\A})$
in non-negative degrees. Suppose that $K_\ast,L_\ast$ are
$\mathbf{G}$-acyclic and $\mathbf G$-presentable. Then any
monoidal natural isomorphism $H_0(K_\ast)\longrightarrow
H_0(L_\ast)$ lifts to a unique isomorphism $f_\ast :K_\ast
\longrightarrow L_\ast$ in $\mathbf{Mon}(\C,
\Chains{\A})[\Sigma^{-1}]$.
\end{corollary}

\subsection{Weak homotopy type} In the following sections we will
apply the general results above to some chain-valued functors
defined in the category of topological spaces or in the category
of simplicial sets. In all cases, the acyclic classes in
$\Chains{\Hom(\C,\A)}$ will come from acyclic classes in
$\Chains{\A}$, and the cotriples in $\Chains{\Hom(\C,\A)}$ will
come from monoidal cotriples in $\C$. Let us describe the main
features of this situation.

\subsubsection{} \label{weakformorphisms}
Let $\Gamma$ be an acyclic class in $\Chains{\A}$ and $\Sigma$ the
associated class of weak equivalences. We extend $\Sigma$ to a
class of weak equivalences $\widetilde{\Sigma}$ in
$\Chains{\Hom(\X,\A)}$ componentwise: a morphism
$f:K_\ast\longrightarrow L_\ast$ is in $\widetilde \Sigma$ if, and
only if, for all objects $X$ of $\C$,
$f_X:K_\ast(X)\longrightarrow L_\ast(X)$ is in $\Sigma$.

If we take the contractible complexes in $\Chains{\A}$ as the
acyclic class $\Gamma$, we will say that $\widetilde{\Gamma}$ is
the acyclic class of {\em weakly contractible} functors in
$\Chains{\Hom(\X,\A)}$. The morphisms in $\widetilde\Sigma$ will
be called {\em weak homotopy equivalences}. If $\Gamma$ is the
class of acyclic complexes, we will say that $\Sigma$ is the class
of {\em weak quasi-isomorphisms}.

\subsubsection{Monoidal cotriples}\label{inducedcotriple}
Given a monoidal category $\C$, a {\em monoidal cotriple} in $\C$
is a cotriple $\mathbf{G}=(G,\varepsilon, \delta)$, such that
$G:\C\longrightarrow \C$ is a monoidal functor and $\varepsilon:
G\Rightarrow \id_\C$ and $\delta: G\Rightarrow G2$ are monoidal
natural transformations.

If $\mathbf{G}=(G,\varepsilon, \delta)$ is a monoidal cotriple in
$\C$, define a functor
$$
\widetilde{G} : \mathbf{Mon}(\C, \Chains{\A}) \longrightarrow
\mathbf{Mon}(\C, \Chains{\A}),
$$
by composition $\widetilde{G}(F_\ast)= F_\ast G$, and let
$\widetilde\varepsilon:\widetilde{G}\Rightarrow \id$ and
$\widetilde\delta : \widetilde{G}\Rightarrow \widetilde{G}^2$ be
the natural transformations induced by $\varepsilon, \delta$. It
follows easily from the definitions that $\widetilde{\mathbf{G}} =
(\widetilde G, \widetilde\varepsilon, \widetilde\delta)$ is a
cotriple in $\mathbf{Mon}(\C, \Chains{\A})$.

\begin{proposition}\label{barra.de.quis}
Let $\Gamma$ be an acyclic class in $\Chains{\A}$ and $\mathbf G$
a monoidal cotriple in $\C$. Consider the class of morphisms
$\widetilde \Sigma$ and the cotriple $\widetilde{\mathbf G}$
induced in $\mathbf{Mon}(\C, \Chains{\A})$. Then,
$\widetilde{\mathbf G}$ is compatible with $\widetilde\Sigma$,
that is, $\widetilde G(\widetilde\Sigma )\subseteq \widetilde
\Sigma$.
\end{proposition}
\begin{proof}
If $f : K_\ast \longrightarrow L_\ast$ is in $\widetilde \Sigma$,
then $f_X: K_\ast (X) \longrightarrow L_\ast (X)$, for all objects
$X$ of $\C$. In particular, $f_{G(X)}\in \Sigma$, thus $\widetilde
G(f)\in \widetilde\Sigma$.
\end{proof}

\subsection{}\label{descenso}
A detailed analysis of the proof of the acyclic models theorem
shows that, in fact, the main ingredients of our proof are a
monoidal functor $s_{EM}:\Delta^o \Chains{\A}\longrightarrow
\Chains{\A}$ and a saturated class of morphisms $\Sigma$ of
$\Chains{\A}$, such that
\begin{itemize}
\item[(1)] $\Sigma$ contains the homotopy equivalences. \item[(2)]
If $f_{\ast \ast} :C_{\ast \ast}\longrightarrow D_{\ast \ast}$ is
a morphism of double complexes such that $f_{n\ast}$ is in
$\Sigma$, then $s_{EM}(f_{\ast \ast})$ is also in $\Sigma$.
\end{itemize}
We will analyze this more general situation elsewhere.

\section{Application: comparison of singular and cubical chains}
In this section we compare the singular and cubical chains of a
topological space in the monoidal setting by applying the results
above, extending the well-known classic comparison theorem,
\cite{EM1}.

\subsection{}\label{seccion43}
Let $\Top$ denote the category of topological spaces, which is a
monoidal category under cartesian product, and let
$\Chains{\mathbb{Z}}$ be the category of chain complexes of
abelian groups. We fix the class $\Sigma$ of weak homotopy
equivalences in the functor category
$\Hom(\Top,\Chains{\mathbb{Z}})$, \ref{weakformorphisms}.

It is well known, \cite{EM1}, that the functor of {\em singular
chains}
$$
S_\ast : \Top\longrightarrow \Chains{\mathbb{Z}},
$$
together with the shuffle product as the K{\"u}nneth morphism is a
monoidal functor.

We can also consider the functor of {\em cubical chains},
$$
C_\ast : \Top\longrightarrow \Chains{\mathbb{Z}}.
$$
For this we take Massey's notations (\cite{Mas}): for a
topological space $X$, let $Q_p(X)$ be the free abelian group
generated by the continuous maps $I^p\longrightarrow X$, where $I$
is the unit interval of the real line. Its elements are the
cubical chains of $X$. Now, $C_n(X)$ is defined as the quotient of
$Q_n(X)$ modulo the degenerated chains. Together with the cross
product
$$
\times : C_\ast(X)\otimes C_\ast(Y)\longrightarrow C_\ast (X\times
Y),
$$
which for singular cubes $c : I^p \longrightarrow X$ and $d : I^q
\longrightarrow Y$ is defined as the cartesian product
$$
c \times d : I^{p+q} = I^p \times I^q \longrightarrow X \times Y,
$$
$C_\ast$ is a monoidal functor.

\subsection{A monoidal cotriple in $\Top$}\label{cotop}

In order to compare $S_\ast$ and $C_\ast$ as monoidal functors
following corollary \ref{unicidadgeneral}, we introduce a monoidal
cotriple in $\Top$. It is a model induced cotriple, (see
\cite{B02}, \S 4.2).

In the classical version of the acyclic models theorem the models
are the standard simplexes $\Delta^m$, $m\geq 0$. To obtain a
monoidal version we take products of these spaces.

Therefore, for any sequence $\underline{n} = (n_1,\dots ,n_r)$,
$n_i\geq 0$, take $\Delta^{\underline{n}} =
\Delta^{n_1}\times\dots\times \Delta^{n_r}$. We define ${\mathbf
G}=(G, \varepsilon ,\delta )$ as the model induced cotriple with
models $\Delta^{\underline n}$. That is, the functor
$G:\Top\longrightarrow \Top$ takes a topological space $X$ to
$$
G(X) = \bigsqcup_{\alpha:\Delta^{\underline{n}}\rightarrow X}
(\Delta^{\underline{n}},\alpha),
$$
where $(\Delta^{\underline{n}},\alpha)$ is a copy of
$\Delta^{\underline{n}}$ indexed by the continuous map $\alpha
:\Delta^{\underline{n}}\rightarrow X$ and where the disjoint union
is over all sequences $\underline{n}$ and maps $\alpha$. $G$ takes
a continuous map $f:X\longrightarrow Y$, to the map
$$
G(f) :\bigsqcup_{\alpha:\Delta^{\underline{n}}\rightarrow X}
(\Delta^{\underline{n}},\alpha) \longrightarrow
\bigsqcup_{\beta:\Delta^{\underline{n}}\rightarrow Y}
(\Delta^{\underline{n}},\beta),
$$
which is the identity from $(\Delta^{\underline{n}},\alpha)$ to
$(\Delta^{\underline{n}},f\circ\alpha)$.

For a topological space $X$ define the map
$$
\varepsilon_X : G(X) =
\bigsqcup_{\alpha:\Delta^{\underline{n}}\rightarrow X}
(\Delta^{\underline{n}},\alpha ) \longrightarrow X,
$$
which over $(\Delta^{\underline{n}},\alpha)$ is
$\alpha:\Delta^{\underline{n}}\longrightarrow X$, so we have a
natural transformation $\varepsilon : G\Rightarrow {\id}$. Finally
we define a natural transformation $\delta:G\Rightarrow G2$: the
iteration of $G$ gives the functor
$$
G2(X) =
\bigsqcup_{\Delta^{\underline{m}}\stackrel{\beta}{\longrightarrow}
\Delta^{\underline{n}}\stackrel{\alpha}{\longrightarrow} X}
((\Delta^{\underline{m}},\alpha),\beta),
$$
so we can define $\delta : G(X)\longrightarrow G2(X)$ as the
identity from $(\Delta^{\underline{n}},\alpha)$ to
$((\Delta^{\underline{n}},\alpha),{\id})$.

Moreover, this cotriple has a monoidal structure. To define the
K\"unneth morphisms of $G$ consider topological spaces $X,Y$ and
take the map
$$
\kappa_{X,Y}: G(X)\times G(Y) =
\left(\bigsqcup_{\alpha:\Delta^{\underline{n}}\rightarrow X}
(\Delta^{\underline{n}},\alpha)\right)\times
\left(\bigsqcup_{\beta:\Delta^{\underline{m}}\rightarrow Y}
(\Delta^{\underline{m}},\beta)\right) \longrightarrow G(X\times Y)
= \bigsqcup_{\gamma:\Delta^{\underline{r}}\rightarrow X\times Y}
(\Delta^{\underline{r}},\gamma),
$$
given by
$$
(\Delta^{\underline{n}},\alpha)\times
(\Delta^{\underline{m}},\beta)\stackrel{\id}{\longrightarrow}
(\Delta^{\underline{n}}\times\Delta^{\underline{m}},\alpha\times\beta).
$$
It is straightforward to prove that $\mathbf G$ is a monoidal
cotriple in $\Top$.

\subsection{} As in previous sections, we denote also by $\mathbf{G}=(G,\varepsilon,
\delta)$ the cotriple induced in
$\mathbf{Mon}(\Top,\Chains{\mathbb{Z}})$ by the monoidal cotriple
in $\Top$ defined above.

\begin{theorem}\label{unicidadtop}
The singular and cubical chain monoidal functors $S_\ast,C_\ast
:\Top\longrightarrow \Chains{\mathbb{Z}}$ are weakly homotopy
equivalent monoidal functors, that is, they are weakly equivalent
in $\mathbf{Mon}(\Top, \Chains{\mathbb Z})$.
\end{theorem}
\begin{proof}
Both $H_0S_\ast$ and $H_0C_\ast$ are the functors that associate
to a topological space $X$ the free group generated by the points
of $X$, so the result will follow from corollary
\ref{unicidadgeneral} after we prove that $S_\ast , C_\ast$ are
$\mathbf G$-presentable and $\mathbf{G}$-acyclic as objects of
$\mathbf{Mon}(\Top,\Chains{\mathbb{Z}})$.

With respect to $\mathbf G$-acyclicity for $S_\ast$ we have to
prove that $G(S_\ast)\longrightarrow G(H_0)$ is a weak
equivalence. That is to say, that the morphisms
$$
\bigoplus_{\alpha:\Delta^{\underline n} \longrightarrow
X}S_\ast(\Delta^{\underline n})\longrightarrow
\bigoplus_{\alpha:\Delta^{\underline n}\longrightarrow X}
\mathbb{Z},
$$
where in each summand the morphism $S_\ast(\Delta^{\underline
n})\longrightarrow \mathbb{Z}$ is the natural augmentation, are
homotopy equivalences, for all $X$. As $\Delta^{\underline{n}}$
are contractible spaces, the result is clear. Similarly, we can
prove the $\mathbf G$-acyclicity of $C_\ast$.

To prove $\mathbf G$-presentability it is sufficient, by
proposition \ref{split}, to prove that $\varepsilon_{C_n}$ and
$\varepsilon_{S_n}$ split, for all $n\geq 0$. We can define a
natural transformation
$$
\theta_S : S_n \longrightarrow G(S_n) ,
$$
by sending a singular simplex $\sigma : \Delta^n\longrightarrow X$
on the topological space $X$ to the simplex of $G(X)$ given by
${\id}:\Delta^n\longrightarrow (\Delta^n,\sigma)$. It follows from
the definition that $\varepsilon\circ \theta_S = \id$, so this
morphism splits the standard resolution $B_\ast(S_\ast)$.

For the cubical chains, observe that $I^n$ is homeomorphic to
$\Delta1\times\stackrel{n)}{\dots}\times \Delta1$, so that a
$n$-cube on a space $X$ is given by a map
$\Delta1\times\stackrel{n)}{\dots}\times \Delta1\longrightarrow
X$. We can now follow the definition of $\theta_S$ to define a
natural transformation for cubical chains
$$
\theta_Q : Q_n \longrightarrow G(Q_n ),
$$
by sending a cubical simplex $c:
\Delta1\times\stackrel{n)}{\dots}\times \Delta1\longrightarrow X$
to the cubical simplex of $G(X)$ given by ${\id}:
\Delta1\times\stackrel{n)}{\dots}\times \Delta1\longrightarrow
(\Delta1\times\stackrel{n)}{\dots}\times \Delta1\longrightarrow X,
c)$. This natural transformation $\theta_Q$ splits the
augmentation $\varepsilon_{Q_n}$. But the natural projection
$\pi_n :Q_n\longrightarrow C_n$ admits a section $\nu_n$,
$\pi_n\nu_n=\id$, (see \cite{Mas}, lemma 5.5), so we can define
$\theta_C = \theta_Q\nu_n$, which splits $\varepsilon_{C_n}$, as
is easily verified.
\end{proof}

\section{Symmetric monoidal functors}

In this section we indicate how to extend the results of the
previous sections to symmetric monoidal functors. We will focus
our attention on the symmetric version of corollary
\ref{unicidadgeneral}, which we will apply to obtain a comparison
result in the symmetric setting for the singular and ordered
cubical functors acting on topological spaces.

\subsection{Acyclic models for symmetric monoidal functors}
Let $\C$ be a symmetric monoidal category, $\A$ a symmetric
abelian monoidal category, $\Gamma$ an acyclic class in
$\Chains{\Hom(\C, \A)}$ and $\Sigma$ the class of weak
equivalences. $\Sigma$ determines a class of morphisms in
$\Symon(\C, \Chains {\A})$, which will be represented by the same
symbol. Let $\mathbf G$ be a cotriple in $\Symon(\C, \Chains
{\A})$ compatible with $\Sigma$, i.e., $G(\Sigma)\subseteq
\Sigma$.

Following the general procedure in \ref{bar}, given $\mathbf{G}$
in $\Symon(\C,\Chains{\A})$, we can associate a functor
$$
B_\bullet: \Symon(\C,\Chains{\A})\longrightarrow
\Delta^o\Symon(\C,\Chains{\A}),
$$
in such a way that the simplicial object $B_\bullet(K_\ast)$ is
augmented to $K_\ast$ by $\varepsilon$. By composing this functor
with the symmetric simple functor, \ref{symsimple}, we obtain a
functor
$$
B_\ast :
\Symon(\C,\Chains{\A})\longrightarrow\Symon(\C,\Chains{\A}).
$$
The natural transformation $\varepsilon$ gives a natural
transformation $B_\ast\Rightarrow \id$.

We can now reproduce the definitions and results from section $\S
4$ in the symmetric setting. In particular, we have the following
result, which is analogous to corollary \ref{unicidadgeneral}.

\begin{theorem}\label{unicidadsim}
Let $K_\ast,L_\ast$ be objects of $\Symon(\C, \Chains{\A})$ in
non-negative degrees. Suppose that $K_\ast,L_\ast$ are
$\mathbf{G}$-acyclic and $\mathbf G$-presentable. Then any
symmetric monoidal natural isomorphism $H_0(K_\ast)\longrightarrow
H_0(L_\ast)$ lifts to a unique isomorphism $f_\ast :K_\ast
\longrightarrow L_\ast$ in $\Symon(\C, \Chains{\A})[\Sigma^{-1}]$.
\end{theorem}

A symmetric monoidal cotriple $\mathbf{G}=(G,\varepsilon, \delta)$
on a symmetric monoidal category $\C$ is a monoidal cotriple such
that $G$ is a symmetric monoidal functor and $\varepsilon, \delta$
are monoidal transformations. As in \ref{inducedcotriple}, if $\D$
is another symmetric monoidal category, a symmetric monoidal
cotriple on $\C$ induces a cotriple on $\Symon(\C,\D)$, that will
also be denoted by $\mathbf{G}$, which satisfies
$G(\Sigma)\subseteq \Sigma$ (see proposition \ref{barra.de.quis}).

\subsection{The Kleisli cotriple}\label{symcotop}
The cotriple in $\Top$ defined in \ref{cotop} is not symmetric.
Therefore in order to apply theorem \ref{unicidadsim} with
$\C=\Top$ we must first define a suitable symmetric monoidal
cotriple on $\Top$. The following cotriple, introduced by Kleisli
in \cite{K2}, will turn out to be symmetric and monoidal: for a
topological space $X$ and an element $x\in X$, let $P(X,x)$ denote
the space of pointed continuous paths $\alpha :
(I,0)\longrightarrow (X,x)$, which is topologized by the
compact-open topology. Then, define $G$ on $X$ by
$$
G(X) = \bigsqcup_{x\in X} P(X,x),
$$
so $G(X)$ is the set of paths in $X$, but with a topology that is
not the path space topology. Nevertheless, note that if $Z$ is a
connected space, a map from $Z$ to $G(X)$ is equivalent to a map
$\beta: Z\times I\longrightarrow X$ such that $\beta(z,0)=x$, for
some $x\in X$ and all $z\in Z$. The action of $G$ on a continuous
map $f:X\longrightarrow Y$ is given by composition, that is, if
$\alpha\in P(X,x)$, then $G(f)(\alpha)=f\circ\alpha $.

Define a natural transformation $\varepsilon : G\Rightarrow \id$
by evaluating paths at $1$. Finally, note that the iteration $G2$
over a space $X$ is given by
$$
G2(X) = \bigsqcup_{x\in X} \bigsqcup_{\alpha\in P(X,x)}
P(P(X,x),\alpha).
$$
An element of $P(P(X,x),\alpha)$ is determined by a map
$\widetilde{\alpha}: I2\longrightarrow X$ with
$\widetilde{\alpha}(t,0)=\alpha(t)$, for all $t\in I$,  and
$\widetilde{\alpha}(0,s) = x$, for all $s\in I$. We define the
natural transformation $\delta : G\Rightarrow G2$ by maps
$P(X,x)\longrightarrow P(P(X,x),x)$ that send $\alpha$ to the map
$\widetilde{\alpha}(t,s)=\alpha(ts)$.

The cotriple $\mathbf{G}=(G,\varepsilon, \delta)$ will be called
the Kleisli cotriple. It is a symmetric monoidal: for topological
spaces $X,Y$ define the K{\"u}nneth morphism
$$
\kappa_{X,Y} : G(X)\times G(Y)=\bigsqcup_{x\in X} P(X,x)\times
\bigsqcup_{y\in Y} P(Y,y) \longrightarrow G(X\times Y)
=\bigsqcup_{(x,y)\in X\times Y} P(X\times Y,(x,y)),
$$
by sending the paths $\alpha\in P(X,x)$ and $\beta\in P(Y,y)$ to
the path $\alpha \times \beta \in P(X\times Y, (x,y))$. It is
evident that these morphisms are compatible with the symmetric
structure of $\Top$.

\subsection{Ordered cubical chains} The cubical chain functor is
not symmetric. However, the ordered cubical chains define a
symmetric monoidal functor that may be compared directly with
$S_\ast$.

\subsubsection{} The ordered cubical chains of a topological space $X$
are defined as follows (see \cite{Ko}). If $c:I^n\longrightarrow
X$ is a singular $n$-cube and $\pi\in\Sigma_n$ define the chain
$\pi c$ as
$$
(\pi c)(t_1,\dots,t_n) = \varepsilon(\pi)\ c(t_{\pi (0)},\dots,
t_{\pi (n)}),
$$
and extend this action to $C_n(X)$ linearly. Let $D_n(X)$ be the
subgroup of $C_n(X)$ generated by chains of the form $c-\pi c$.
Then $D_\ast(X)$ is a subcomplex of $C_\ast(X)$, so we can define
the {\em ordered cubical chains} of $X$ as the quotient complex
$$C^{ord}_\ast(X) = C_\ast(X)/D_\ast(X).$$
There is a natural transformation $C_\ast\Rightarrow
C^{ord}_\ast$.

\subsubsection{} The monoidal structure of the cubical chain
functor $C_\ast$ with the usual cross product is carried over the
quotient by $D_\ast$. In fact, if $X,Y$ are topological spaces,
the cross product maps $C_\ast(X) \otimes D_\ast(Y)$ to
$D_\ast(X\times Y)$ since if $c\in C_p(X), d\in C_q(Y)$ are
singular cubes and $\pi\in\Sigma_p$, then
$$
(c-\pi c) \times d = (c\times d) - \widetilde{\pi}{(c\times d)},
$$
where $\widetilde \pi$ is the element of $\Sigma_{p+q}$ that acts
as $\pi$ on the first $p$ elements and fixes the rest.
Analogously, $\times$  maps $D_\ast(X) \otimes C_\ast(Y)$ to
$D_\ast(X\times Y)$.

\begin{proposition}
The functor $C_\ast^{ord}:\Top\longrightarrow \Chains{\mathbb{Z}}$
is a symmetric monoidal functor and the natural transformation
$C_\ast\Rightarrow C^{ord}_\ast$ is a monoidal natural
transformation.
\end{proposition}
\begin{proof}
The proof follows immediately from the definitions. Note that the
cross product in $C_\ast^{ord}$ is symmetric, since if $X,Y$ are
spaces and $c\otimes d\in C_p^{ord}(X)\otimes C_q^{ord}(X)$, then
\begin{eqnarray*}
\tau(c\otimes d) (t_1,\dots,t_{p+q}) &=& (-1)^{pq}(d\otimes c)
(t_1,\dots,t_{p+q})\\
&=& (d\otimes c) (t_q,\dots, t_{p+q}, t_1,\dots ,t_p)\\
&=& \tau_{X,Y}(c\times d)(t_1,\dots,t_{p+q}),
\end{eqnarray*}
where in the second equality we used the invariance of the
oriented cubical chains by the action of the symmetric group.
\end{proof}

\subsubsection{} We fix in $\Symon(\Top,
\Chains{\mathbb{Z}})$ the class $\Sigma$ of weak homotopy
equivalences, that is, the symmetric monoidal functors $f$ such
that $f_X$ is a homotopy equivalence for each space $X$. Let
$\mathbf G$ be the cotriple induced by the Kleisli cotriple, which
satisfies $G(\Sigma)\subseteq \Sigma$, see \ref{barra.de.quis}.

\begin{theorem}\label{unicidadtopsim}
The singular and ordered cubical chain functors $S_\ast,
C_\ast^{ord} : \Top\longrightarrow \Chains{\mathbb{Z}}$ are weakly
homotopy equivalent symmetric monoidal functors.
\end{theorem}
\begin{proof}
For a topological space $X$, both groups $H_0(S_\ast(X))$ and
$H_0(C_\ast^{ord}(X))$ are isomorphic to the free group generated
by the points of $X$, so the result will follow from theorem
\ref{unicidadsim} after we prove that $S_\ast$ and $C_\ast^{ord}$
are $\mathbf G$-presentable and $\mathbf G$-acyclic with respect
to the Kleisli cotriple. This has been proved by Barr (see
\cite{B02}) for the singular chains functor. Let us prove it for
the ordered cubical chains.

Firstly we prove the $\mathbf G$-acyclicity of $C_\ast^{ord}$. We
will prove that there is a chain contraction $s$ for the complex
$C_\ast G\longrightarrow H_0G$ such that for any singular $n$-cube
$c$ and any $\sigma\in\Sigma_n$, $s(\sigma c) = \sigma s(c)$.

Let $X$ be any topological space and $c:I^n\longrightarrow G(X)$ a
singular $n$-cube. By the connectedness of the standard cube,
there is a point $x\in X$ such that the map $c$ factors through a
map $c:I^n\longrightarrow P(X,x)$. By adjunction, $c$ is
equivalent to a map $\widetilde{c} : I^{n+1}\longrightarrow X$,
which satisfies $\widetilde{c}(t_1,\dots ,t_n ,0)=x$, for
$(t_1,\dots , t_n)\in I^n$. Taking into account the product
decomposition $I^{n+1}=I^n\times I$, we will write the value of
$\widetilde c$ at the point $(t_1,\dots, t_n)\in I^n$ and $u\in I$
by $\widetilde{c}(t_1,\dots,t_n;u)$.

Note that $H_0(G(X))$ is the free group generated by the elements
of $X$, since the spaces $P(X,x)$ are contractible, and define
$s:H_0(G(X))\longrightarrow C_0(G(X))$ by $s(x)=p_x$, where $p_x$
denotes the constant path at $x$.

For $n\geq 0$ define $\mu:I^{n+1}\times I\longrightarrow I^n\times
I$ by
$$
\mu(t_1,\dots,t_{n+1};u)=(t_1,\dots,t_n;t_{n+1}u)
$$
and  define $s:Q_n(G(X))\longrightarrow Q_{n+1}(G(X))$ by
$$
s(\widetilde c) = (-1)^{n+1}\widetilde{c}\circ\mu.
$$
Next, we recall the definition of the differential $d:Q_n(X)
\longrightarrow Q_{n-1}(X)$ of the cubical chain complex: for
$1\le i\le n,\ \epsilon\in\{0,1\}$, let
$\delta_i^{\epsilon}:I^{n-1} \longrightarrow I^{n}$ denote the
face defined by $ \delta^{\epsilon}_{i}
(t_1,\dots,t_{n-1})=(t_1,\dots,\epsilon,\dots,t_{n-1})\,\,, $
where $\epsilon$ is in the $i$-th place; if $c\in Q_n(X)$, $d(c)$
is defined by
$$
d (c)=\sum_{i,\epsilon}(-1)^{i+\epsilon}c\circ
\delta_i^{\epsilon}.
$$
It is clear that on $Q_{n+1}$ we have
$\delta_i^\varepsilon\circ\mu = \mu\circ\delta_i^\varepsilon$, for
$1\leq i\leq n$, and that
\begin{eqnarray*}
(\mu\circ\delta_{n+1}^0)(t_1,\dots,t_n;u) &=& (t_1,\dots,t_n,0),\\
(\mu\circ\delta_{n+1}^1)(t_1,\dots,t_n;u) &=& (t_1,\dots,t_n;u),
\end{eqnarray*}
thus, an easy calculation proves that
$$
(ds+sd)(\widetilde{c}) = \widetilde c - x.
$$
As $s$ sends degenerated chains to degenerated chains, it defines
a map $s:C_n(G(X))\longrightarrow C_{n+1}(G(X))$, and as $x$
represents the constant path $p_x$, which is degenerated, the
equality above reduces to $ds+sd=\id$ in $C_\ast(G(X))$.

The contraction $s$ is trivially compatible with the action of the
symmetric group, so it defines a contraction for $C_\ast^{ord}$,
and as a consequence $C_\ast^{ord}$ is $\mathbf G$-acyclic.

To prove $\mathbf G$-presentability, notice that the natural
transformation $\theta_n : Q_n\longrightarrow Q_nG$ given by
$$
\theta_n(\sigma)(t_1,\dots ,t_n)(s) = c(st_1,\dots, st_n),
$$
is a section of $\varepsilon_Q $ which is compatible with taking
quotients modulo degenerated chains and the action of the
symmetric group $\Sigma_n$, so it defines a section
$\theta_n^{ord} : C_n^{ord}\longrightarrow C_n^{ord}G$, and we can
apply proposition \ref{split}.
\end{proof}

\begin{remark} \label{cadenesnormalitzades} Note that if
$S^N_\ast : \Top\longrightarrow \Chains{\mathbb{Z}}$ denotes the
functor of normalized singular chains, which is also a symmetric
monoidal functor since the degenerated singular chains are
invariant by shuffle product, (see \cite{EM2}), then the
projection $S_\ast\Rightarrow S^N_\ast$ is a monoidal natural
transformation of symmetric monoidal functors which is a weak
equivalence. This is a classical result and follows directly from
the fact that the degenerated singular chains are a direct factor
of singular chains.
\end{remark}

\section{Application to operads}

\subsection{Operads} Let us recall some definitions and notations
about operads (see \cite{MSS}).

\subsubsection{} Let $\mathbf \Sigma$ be the {\em symmetric groupoid},
that is, the category whose objects are the sets $\{1,\dots, n\}$,
$n\geq 1$, and the only morphisms are those of the symmetric
groups $\Sigma_n$.

\subsubsection{} Let $\C$ be a monoidal category. The category of
contravariant functors from $\mathbf \Sigma$ to $\C$ is called the
category of $\mathbf\Sigma$-modules and is represented by
$\sMod_{\mathcal{C}}$.We identify its objects with sequences of
objects in $\C$, $E=\left((E(l)\right)_{l\geq 1}$, with a right
$\Sigma_l$-action on each $E(l)$. If $E$ and $F$ are
$\mathbf\Sigma$-modules, a {\it morphism of
$\mathbf\Sigma$-modules} $f:E \longrightarrow F$ is a sequence of
$\Sigma_l$-{ equivariant\/} morphisms $f(l): E(l) \longrightarrow
F(l), \ l\geq 1$.

\subsubsection{} A {\em unital $\mathbf\Sigma$-operad\/} (an {\em
operad\/} for short) in $\mathcal{C}$ is a $\mathbf \Sigma$-module
$ P$ together with a family of {\em composition morphisms\/}
$$\gamma_{l;m_1,\dots ,m_l} : P(l)\otimes P(m_1) \otimes \dots
\otimes P(m_l) \longrightarrow P(m_1 + \cdots + m_l),$$ and a {\em
unit morphism\/} $$\eta : \mathbf{1} \longrightarrow P(1),$$ which
satisfies the axioms of equivariance, associativity and unit. A
{\em morphism of operads\/} is a morphism of
$\mathbf\Sigma$-modules that is compatible with structure
morphisms. We use $\Op_{\mathcal{C}}$ to denote the category of
operads in $\mathcal C$ and its morphisms.

An operad in $\Top$ is called a {\em topological operad}. If $R$
is a ring, an operad in $\Chains{R}$ is called a {\em dg operad}.
We are especially interested in the cases $R=\mathbb{Z},
\mathbb{Q}$.

\subsection{} We can now extend the comparison result between
$S_\ast$ and $C_\ast^{ord}$ to topological operads.

\subsubsection{} If $F:\C\longrightarrow \D$ is a symmetric monoidal
functor between monoidal categories, it is easy to prove that,
applied componentwise, $F$ induces a functor between the
corresponding categories of operads
$$
\Op_F : \Op_\C \longrightarrow \Op_\D ,
$$
also denoted by $F$. Therefore, singular and cubical ordered
chains induce functors
$$
S_\ast, C_\ast^{ord}: \Op_\Top\longrightarrow
\Op_{\Chains{\mathbb{Z}}}.
$$

Moreover, if $\D$ has a notion of weak equivalence, it extends to
$\Op_\D$ componentwise. For instance, we can consider the weak
equivalence relation in the category of chain complexes
$\Chains{R}$ induced by quasi-isomorphisms. Now it easily follows
from the definitions that

\begin{proposition}
Let $\C$ be a monoidal category, $R$ a ring and $F,G\in
{\mathbf{Mon}}(\C,\Chains{R})$. If $F,G$ are weakly equivalent
(with respect to quasi-isomorphism in $\Chains{R}$), then the
functors
$$\Op_F, \Op_G: \Op_\C\longrightarrow \Op_{\Chains{R}}$$ are
weakly equivalent.
\end{proposition}

This proposition together with theorem \ref{unicidadtopsim}
applied to $S_\ast, C_\ast^{ord}$ gives

\begin{theorem}\label{unicidadoperad}
The functors
$$
S_\ast, C_\ast^{ord} : \Op_\Top\longrightarrow
\Op_{\Chains{\mathbb{Z}}}
$$
are weakly equivalent (with respect to quasi-isomorphism). In
particular, for a topological operad $P$ the dg operads
$S_\ast(P), C_\ast^{ord} (P)$ are weakly equivalent.
\end{theorem}

If $\rm{Ho}\Op_{\Chains{\mathbb{Z}}}$ denotes the localization of
$\Op_{\Chains{\mathbb{Z}}}$ with respect to quasi-isomorphisms, we
obtain
\begin{corollary}
There is an isomorphism of functors
$$
S_\ast \cong C_\ast^{ord} : \Op_\Top\longrightarrow \mbox{\rm
Ho}\Op_{\Chains{\mathbb{Z}}}.
$$
\end{corollary}

\subsubsection{} Theorem \ref{unicidadoperad} may be applied to
compare the categories of $S_\ast(P)$ and $C_\ast^{ord}(P)$
algebras up to homotopy.

Recall that, given a dg operad $P$, a $P$-algebra is a chain
complex $V$ of finite type together with a morphism
$P\longrightarrow \E[V]$, where $\E[V]$ denotes the operad of
endomorphisms of $V$. Following \cite{GNPR}, definition 7.3.1,
define a {\em $P$-algebra up to homotopy} as a finite type complex
$V$ together with a morphism $P\longrightarrow \E[V]$ in
$\mbox{\rm Ho}\Op_{\Chains{\mathbb{Z}}}$. From theorem
\ref{unicidadoperad}, it immediately follows that

\begin{corollary}\label{uptohomotopy}
Let $P$ be a topological operad. The categories of $S_\ast(P)$ and
$C^{ord}_\ast(P)$ algebras up to homotopy are equivalent.
\end{corollary}

This remark applies to Deligne's conjecture, which is commonly
expressed as follows: let $A$ be an associative algebra over a
ring $R$, and let $C^\ast (A;A)$ be the associated Hochschild
complex.

{\bf Deligne's conjecture.} {\em For any associative $R$-algebra
$A$, the complex $C^\ast (A;A)$ is naturally an algebra over the
singular chains of the little discs operad $\D_2$ or a suitable
version of it}.

There is some ambiguity in this statement with reference to the
chain model of $\D_2$ to be used and whether obtaining a solution
for one model means that a solution is obtained for any other
model. From our results it follows that there is only one chain
model of a topological operad up to homotopy, therefore, this
ambiguity disappears when we work with structures up to homotopy,
linking in this way the solutions of Deligne's conjecture given by
Kontsevich, who used $C_\ast^{ord}(\D_2)$, (see \cite{Ko}), and
MacClure-Smith, Tamarkin and Voronov, who used singular chains
(see, for example, \cite{T}).

\subsection{Formality} Recall that a dg operad
$P$ is {\em formal} if it is quasi-isomorphic to its homology
operad $H(P)$, (see \cite{MSS}). From theorem \ref{unicidadoperad}
we can deduce the following:

\begin{corollary}\label{formalidad.indistinta}
Let $P$ be a topological operad. Then $S_\ast(P ;\mathbb{Q})$ is a
formal operad if, and only if, $C_\ast^{ord} (P ;\mathbb{Q})$ is a
formal operad.
\end{corollary}

M. Kontsevich proved the formality of the ordered cubical chains
of the little $k$-discs operad after tensoring by the real
numbers, $\mathbb R$, (see \cite{Ko}, theorem 2). The independent
nature of formality on the base field proved in \cite{GNPR}
implies that $C_\ast^{ord}(\D_k; \mathbb{Q})$ is also formal over
the rational field, so by corollary \ref{formalidad.indistinta},
$S_\ast(\D_k;\mathbb{Q})$ is formal.

\subsection{Modular operads}
The results above may be extended to chain models of modular
operads. We refer to \cite{GK} and \cite{MSS} for the definitions
concerning modular operads. Given a symmetric monoidal category we
denote by $\MOp_\C$ the category of modular operads of $\C$.

As for operads, it follows from the definitions that every
symmetric monoidal functor $F : \mathcal{C} \longrightarrow
\mathcal{D}$ applied componentwise induces a functor
$$
\MOp_{F}:\MOp_{\mathcal{C}} \longrightarrow \MOp_{\mathcal{D}}.
$$
In particular, the singular and cubical chain functors of
topological spaces extend to functors defined in the category of
topological modular operads. Now, as in theorem
\ref{unicidadoperad}, one has

\begin{theorem}\label{unicidadoperadmodular}
The functors
$$
S_\ast, C_\ast^{ord} : \MOp_\Top\longrightarrow
\MOp_{\Chains{\mathbb{Z}}}
$$
are weakly equivalent. In particular, for a topological modular
operad $P$, the dg modular operads $S_\ast(P),  C_\ast^{ord} (P)$
are weakly equivalent.
\end{theorem}
From theorem \ref{unicidadoperadmodular} it immediately follows
that
\begin{proposition}
Let $P$ be a topological modular operad. Then $S_\ast(P
;\mathbb{Q})$ is a formal modular operad if and only if
$C_\ast^{ord} (P ;\mathbb{Q})$ is a formal modular operad.
\end{proposition}

We can apply this result to the modular operad $\M$: the family
$\overline{\M}((g,l)) = \overline{\M}_{g,l}$ of
Deligne-Knudsen-Mumford moduli spaces of stable genus $g$
algebraic curves with $l$ marked points, with the maps that
identify marked points, is a modular operad in the category of
projective smooth Deligne-Mumford stacks, \cite{GK}. In
\cite{GNPR} we proved that $S_\ast(\overline{\M};\mathbb{Q})$ is a
formal modular operad, so any other chain model is also formal.
For instance, for cubical chains we can state that

\begin{corollary}
The dg modular operad $C_\ast^{ord}(\overline{\M};\mathbb{Q})$ is
a formal modular operad.
\end{corollary}

\section{Contravariant functors}

If we work in the category of contravariant monoidal functors
between monoidal categories, $\mathbf{Mon}(\C^o,\D)$, a cotriple
on $\C$ induces a triple on $\mathbf{Mon}(\C^o,\D)$, so that the
standard construction produces cosimplicial objects instead of
simplicial objects. As remarked in \ref{comon=mon} this can be
avoided by identifying such functors with covariant comonoidal
functors between $\C$ and $\D^o$. However, for functors with
values in $\Cochains{\mathbb{Z}}$, which are the functors that
appear in the study of the cohomology of topological spaces, the
dual category $\Cochains{\mathbb{Z}}^o$ is not the same as
$\Cochains{\mathbb{Z}}$, so the acyclic models theorem
\ref{unicidadgeneral} has to be appropriately dualized to cover
this case.

In this section we present the minor modifications of the
constructions and results of \S 4 that are necessary to cover this
situation, and we apply them to compare the singular and cubical
cochain functors.

\subsection{The standard construction} Let $\mathbf{T}=(T,\eta ,\mu)$
be a triple (or {\em monad}, see \cite{ML}) in a category $\X$.
There is a functor associated to $\mathbf T$, which we call the
{\em standard construction},
$$
B^\bullet : \X\longrightarrow \Delta \X,
$$
defined by $B^n(X)=T^{n+1}(X)$ and with faces and degeneracies
defined analogously as in \ref{bar}.

\subsubsection{} Let $\C$ be a monoidal category and $\A$ an
additive monoidal category, and $\Sigma$ a class of weak
equivalences of $\mathbf{Mon}(\C^o,\Cochains{\A})$ which contains
the homotopy equivalences and is compatible with the
Alexander-Whitney functor $s_{AW}$. A monoidal cotriple $\mathbf
{G}=(G,\varepsilon, \delta)$ in $\C$ induces a triple $\mathbf T$
on $\X=\mathbf{Mon}(\C^o,\Cochains{\A})$, whose endofunctor is
given by $T(K^\ast) = K^\ast G$, and with $\eta, \mu$ induced by
$\varepsilon, \delta$. By the dual of proposition
\ref{barra.de.quis}, $\mathbf T$ is compatible with $\Sigma$, that
is, $T(\Sigma)\subseteq \Sigma$.

We can compose the standard cosimplicial construction with the
Alexander-Whitney simple functor defined in \ref{AWcontra} to
obtain a functor
$$
B^\ast : \mathbf{Mon}(\C^o,\Cochains{\A}) \longrightarrow
\mathbf{Mon}(\C^o, \Cochains{\A}).
$$
Moreover, the natural transformation $\varepsilon$ of $\mathbf G$
defines a natural transformation $\varepsilon:{\id}\Rightarrow
B^\ast$.

\subsection{Acyclic models for contravariant monoidal functors}
By analogy to section $\S 4$, it is clear how to define the
notions of $\mathbf T$-presentability and $\mathbf T$-equivalence
in $\mathbf{Mon}(\C^o,\Cochains{\A})$, so that we have all the
ingredients to transpose the proof of
\ref{representable=cofibrante} and its corollaries to the category
of contravariant monoidal functors between $\C$ and
$\Cochains{\A}$. In particular, we can state the contravariant
version of corollary \ref{unicidadgeneral} referred to the
assumptions made in the previous paragraph:

\begin{theorem}\label{rep=cof2}
Let $K^\ast, L^\ast$ be objects of
$\mathbf{Mon}(\C^o,\Cochains{\A})$ concentrated in non-negative
degrees, and $\mathbf T$ a triple in
$\mathbf{Mon}(\C^o,\Cochains{\A})$ induced by a monoidal cotriple
$\mathbf G$ on $\C$. Suppose that $K^\ast, L^\ast$ are $\mathbf
T$-acyclic and $\mathbf T$-presentable. Then any monoidal
transformation $H0(K^\ast)\longrightarrow  H0(L^\ast)$ lifts to a
unique morphism $K^\ast\longrightarrow L^\ast$ in
$\mathbf{Mon}(\C^o,\Cochains{\A})[\Sigma^{-1}]$. In particular, if
$H0(K^\ast)$ and $H0(L^\ast)$ are isomorphic, then $K^\ast$ and
$L^\ast$ are weakly homotopy equivalent.
\end{theorem}

\subsection{Application to singular and cubical cochains} For a
topological space $X$, let $S^\ast(X), C^\ast(X)$ denote,
respectively, the complexes of singular and of cubical cochains
defined on $X$:
\begin{eqnarray*}
S^\ast (X) &=& \mbox{Hom} (S_\ast(X),\mathbb{Z}),\\
C^\ast (X) &=& \mbox{Hom} (C_\ast(X),\mathbb{Z}).
\end{eqnarray*}
They define contravariant functors
$$
S^\ast, C^\ast : \Top^o\longrightarrow \Cochains{\mathbb{Z}}.
$$
By dualizing the Alexander-Whitney morphism for $S_\ast$, $S^\ast$
becomes a contravariant monoidal functor. There is also an
explicit Alexander-Whitney associative morphism for cubical chains
$C_\ast(X\times Y)\longrightarrow C_\ast(X)\otimes C_\ast(Y)$ (see
\cite{Mas}, XI\S 5 and exercise XIII.5.1), and as a consequence
$C^\ast$ is also a contravariant monoidal functor.

The monoidal cotriple $\mathbf G$ on $\Top$ defined in \ref{cotop}
induces a triple in $\mathbf{Mon} (\Top^o,\Cochains{\mathbb{Z}})$,
denoted by $\mathbf T$, and it is easily seen (compare
\ref{unicidadtop}) that $S^\ast$ and $C^\ast$ are $\mathbf
{T}$-presentable and $\mathbf T$-acyclic. Thus we are able to
apply theorem \ref{rep=cof2} and deduce:

\begin{theorem}
The singular and cubical cochain functors $S^\ast, C^\ast :
\Top^o\longrightarrow \Cochains{\mathbb{Z}}$ are weakly homotopy
equivalent contravariant monoidal functors, that is, they are
weakly homotopy equivalent in
$\mathbf{Mon}(\Top^o,\Cochains{\mathbb{Z}})$.
\end{theorem}

\section{Application to cohomology theories}
In this section we apply the acyclic models theorem \ref{rep=cof2}
to compare cohomology theories on simplicial sets arising from
simplicial differential graded algebras over a ring $R$.

\subsection{Cohomology theories} Denote by $\Sets$ the category of
simplicial sets. If $A_{\bullet}^{\ast}$ is a simplicial
differential graded $\mathbb Z$-algebra, then to any simplicial
set $X$ we can associate a differential graded algebra $A^\ast(X)$
by taking morphisms, in $\Sets$, from $X$ to $A_\bullet^\ast$. In
order to obtain a good cohomology theory we impose some conditions
on $A_\bullet^\ast$, following Cartan (see \cite{C}).

\subsubsection{} Let $R$ be a ring and $A_{\bullet}^{\ast}$ a
simplicial differential graded $R$-algebra. We will assume that
$A_{\bullet}^{\ast}$ satisfies the following axioms:
\begin{itemize}
\item[(a)] {\em Homology axiom}. For each $p\geq 0$, the natural
morphism $R\longrightarrow A_p^{\ast}$ is a homotopy equivalence
of cochain complexes. In particular, it is a quasi-isomorphism.

\item[(b)] {\em Homotopy axiom}. For each $q\geq 0$, the
simplicial set $A^q_\bullet$ is (simplicially) contractible, i.e.
the homotopy groups of $A^q_\bullet$ are zero.

\item[(c)] {\em Freeness axiom}. $R$ is a principal ideal domain
(PID) and, for all $p,q$, the $R$-module $A_p^q$ is free.
\end{itemize}

\subsubsection{} Associated to $A_\bullet^\ast$ there is a {\em cohomology theory}
defined by
$$
A^\ast (X) = {\rm Hom}_\Delta(X_\bullet ,A^\ast_\bullet),
$$
where $\rm Hom_\Delta$ stands for the homomorphism set in the
simplicial category. $A^\ast(X)$ is a cochain complex, which in
degree $q$ is equal to $\Sets(X_\bullet,A^q_\bullet)$.

\begin{example}
Given a ring $R$, the singular cochain complex $S^\ast(\ ;R)$ is
an example of cohomology theory in the sense above. To see this it
suffices to take the simplicial $R$-algebra $S^\ast_\bullet (R)$
which in simplicial degree $p$ is the cochain $R$-algebra of the
simplicial set represented by $p$, $\Delta[p]$, with evident face
and degeneracies.
\end{example}

Cartan proves in \cite{C}, theorem 1 (see also \cite{Maj}), that
the cohomology of cohomology theories associated with simplicial
differential graded algebras $A^\ast_\bullet$, which satisfies the
homology and homotopy axioms, are isomorphic, and that this
isomorphism comes from a true morphism of complexes if
$A^\ast_\bullet$ satisfies the normalization  $A_p^q = 0$ if
$q>p$, (loc. cit. theorem 2). Moreover, it is compatible with
products if the cohomology theories satisfy some flatness
conditions (loc. cit. theorem 3). We will apply \ref{rep=cof2} to
obtain a comparison result at the chain level, which is stronger
than Cartan's theorems.

Note that in  \cite{Man}, M.Mandell obtains uniqueness results for
cochain theories that satisfy a different set of axioms related to
the classical Eilenberg-Steenrod axioms.

\subsubsection{} The category of simplicial sets $\Sets$
is a monoidal category under the cartesian product. Using the
algebra structure of $A_\bullet^\ast$ it follows that the
cohomology theory $A^\ast$ defines a contravariant monoidal
functor
$$
A^\ast: (\Sets)^o\longrightarrow \Cochains{R}.
$$
\begin{theorem}\label{unicidadderham}
Let $A^\ast_\bullet, B^\ast_\bullet$ be simplicial differential
graded $R$-algebras satisfying axioms (a)-(c). Then the
contravariant monoidal functors
$$
A^\ast , B^\ast : (\Sets )^o\longrightarrow \Cochains{R},
$$
are weakly quasi-isomorphic in the category $\mathbf{Mon}((\Sets
)^o,\Cochains{R})$. In particular, any such cohomology theory is
weakly quasi-isomorphic to the singular cochain complex functor
$S^\ast(\ ;R)$.
\end{theorem}
\begin{proof}
Let $\mathbf G$ be the cotriple in $\Sets$ which is defined, as in
the topological case, by
$$
G(X) =\sqcup_\alpha (\Delta[\underline{n}], \alpha),
$$
where $\Delta[\underline{n}]=\Delta[n_1]\times\dots\times
\Delta[n_r]$. $\mathbf G$ induces a triple
$\mathbf{T}=(T,\eta,\mu)$ in $\mathbf{Mon} (\Sets^o,
\Cochains{R})$.

By hypothesis, $H0(A^\ast(X))$ and $H0(B^\ast(X))$ are isomorphic
to $ H0(X)$, so in order to apply the acyclic models theorem
\ref{rep=cof2} we have to prove that any cohomology theory
$A^\ast$ that satisfies axioms (a)-(c) is $\mathbf T$-acyclic and
$\mathbf T$-presentable, the case of $S^\ast(\ ;R)$ being well
known. The $\mathbf T$-acyclicity will follow from the acyclicity
of $A_\bullet^\ast$ with respect to the differential degree, while
the $\mathbf T$-presentability will follow from the
contractibility with respect to the simplicial degree.

Let us first prove that $A^\ast$ is $\mathbf T$-acyclic. We have
to prove that for any sequence $\underline n = (n_1,\dots ,n_r)$
the complex
$$
\dots \longrightarrow A^q(\Delta^{\underline{n}})\longrightarrow
A^{q-1}(\Delta^{\underline{n}})\longrightarrow \dots
\longrightarrow A0(\Delta^{\underline{n}})\longrightarrow
H0(A^\ast(\Delta^{\underline{n}})) \longrightarrow 0,
$$
is acyclic. We prove this statement for any contractible
simplicial set $X$, $X=\Delta^{\underline{n}}$ being a special
case.

Note that the $q$-cocycles of $A^\ast(X)$ are equal to ${\rm Hom}
(X,Z^qA)$, where $Z^qA$ denotes the simplicial group of
$q$-cocycles of $A_\bullet^\ast$. By the homology axiom
$d:A_\bullet^{q-1}\longrightarrow Z^{q}A$ is surjective, hence it
is a Kan fibration with fiber $Z^{q-1}A$, and it follows also that
$Z^qA$ is connected.

Take a $q$-cocycle $f: X\longrightarrow Z^{q}A$. As $X$ is
contractible, $f$ is homotopic to zero and $Z^qA$ is connected,
$f$ admits an extension to a morphism $X\longrightarrow
A_\bullet^{q-1}$, and consequently it is a boundary.

Let us now turn to $\mathbf T$-presentability. By the
contravariant version of proposition \ref{split}, if each $A^q$,
$q\geq 0$, is $\mathbf T$-split then $A^\ast$ will be $\mathbf
T$-presentable. Therefore we want to define, for each $q$, a
natural transformation $\theta: T(A^q)=A^qG\Rightarrow A^q$ such
that $\theta\eta={\id}$, where for a simplicial set $X$,
$$
\eta_X: A^q(X)\longrightarrow A^qG(X)= \prod_{\underline{n}}
\prod_{\alpha:\Delta^{\underline{n}}\rightarrow
X}A^q(\Delta^{\underline{n}}, \alpha)
$$
is the morphism given by composition: $\eta_X(w) = (w\circ \alpha,
\alpha)$.

We denote by $X(p)$ the $p$-skeleton of the simplicial set $X$, so
$X=\varinjlim (X(p))$. First of all, note that the functors  $A^q$
and $A^qG$ are compatible with the skeleton decomposition, that
is, we have
\begin{eqnarray*}
A^q(X) &=& \varprojlim A^q(X(p)),\\
A^qG(X) &=& \varprojlim A^qG(X(p)).
\end{eqnarray*}
The first isomorphism results from the compatibility of Hom
functors with limits, while for the second we observe that $G$
commutes with filtered colimits.

Thus, to define $\theta$ it is sufficient to define morphisms
$\theta_p : A^qG(X(p))\longrightarrow A^q(X(p))$, $p\geq 0$, which
are sections of $\eta_{X(p)}$, and such that the diagrams
$$
\begin{CD}
A^qG(X(p))
@>\theta_p >> A^q(X(p))\\
@VVV  @VVV\\
A^qG(X(p-1)) @>{\theta_{p-1}}>> {A^q(X(p-1)),}
\end{CD}
$$
where the vertical morphisms are induced by the inclusion
$X(p-1)\subseteq X(p)$, are commutative.

We will define $\theta_p$ inductively on $p$. The inductive step
will be based on the following auxiliary result:

{\bf Extension lemma.} {\em For any $p$, the inclusion
$\partial\Delta[p]\longrightarrow \Delta[p]$ induces a surjection
$f_p:A^q(\Delta[p])\longrightarrow A^q (\partial\Delta[p])$ which
has a $R$-linear section $s: A^q(\partial\Delta[p])\longrightarrow
A^q (\Delta[p])$.}

In fact, note that, as a consequence of the homotopy axiom, for
any simplicial subset $Y\subseteq X$ the induced morphism
$A^q(X)\longrightarrow A^q(Y)$ is surjective, since $A^q_\ast$ is

a contractible Kan complex and, as the inclusion $Y\subseteq X$ is
a cofibration, any map $Y\longrightarrow A^q_\ast$ extends to a
map $X\longrightarrow A^q_\ast$, (see \cite{BG}, and \cite{FHT}
for a more elementary proof).

In particular, the morphism $f_p:A^q(\Delta[p])\longrightarrow A^q
(\partial\Delta[p])$ is surjective. But $A^q (\partial\Delta[p])$
is $R$-projective, because $R$ is a PID and $A^q
(\partial\Delta[p])$ is a submodule of the free $R$-module
$(A^q_{p-1})^{p+1}$. Thus, $f_p$ has an $R$-linear section $s:
A^q(\partial\Delta[p])\longrightarrow A^q (\Delta[p])$.

Let us now return to the inductive definition of $\theta_p$. For
$p=-1$ there is nothing to prove, so let us assume that
$\theta_{p-1}$ has been constructed. The $p$-skeleton $X(p)$ is
obtained from the $(p-1)$-skeleton $X(p-1)$ by the pushout diagram
$$
\begin{CD}
\sqcup \partial\Delta[p] @>>> X(p-1)\\
@VVV @VVV\\
\sqcup\Delta[p]@>>> X(p),
\end{CD}
$$
where the disjoint union is over all nondegenerated maps
$\Delta[p]\longrightarrow X(p)$. Since $A^q$ transforms pushouts
to pullbacks, we get the pullback diagram
$$
\begin{CD}
A^q(X(p))@>>> \prod A^q(\Delta[p])\\
@VVV @VVV \\
A^q(X(p-1)) @>>>\prod A^q(\partial\Delta[p]).
\end{CD}
$$
Analogously, by applying $A^qG$ we obtain the pullback diagram
$$
\begin{CD}
A^qG(X(p))@>>> \prod A^qG(\Delta[p])\\
@VVV @VVV\\
A^qG(X(p-1)) @>>>\prod A^qG(\partial\Delta[p]).
\end{CD}
$$
By induction we have morphisms $\theta_{p-1}=\theta_{X(p-1)}$ and
$\theta'_{p-1}=\theta_{\partial\Delta[p]}$, so it suffices to
define $\theta'_p=\theta_{\Delta[p]}$ to make the diagram
$$
\begin{CD}
A^qG(\Delta[p]) @>g_p>> A^qG(\partial\Delta[p]) @<<< A^qG(X(p-1))\\
@V\theta'_pVV @V\theta'_{p-1}VV @V{\theta_{p-1}}VV\\
A^q(\Delta[p]) @>f_p>> A^q(\partial\Delta[p]) @<<< A^q(X(p-1))
\end{CD}
$$
commutative and which is a section of $\eta_{\Delta[p]}$, (the
square on the right commutes by induction).

Let $\pi: A^qG(\Delta[p])\longrightarrow A^q(\Delta[p])$ be the
natural projection morphism, so that $\pi\eta =\id$, and let $s$
be a section of $f_p$, $f_ps=\id$, whose existence is guaranteed
by the extension lemma above. We define the morphism
$$
\theta'_p = \pi + s(\theta'_{p-1}g_p - f_p\pi).
$$
Commutativity easily follows, since
$$
f_p\theta'_p = f_p\pi + f_ps(\theta'_{p-1}g_p - f_p\pi) = f_p\pi +
\theta'_{p-1}g_p - f_p\pi = \theta'_{p-1}g_p.
$$
Moreover, $\theta'_p$ is a section of $\eta_{\Delta[p]}$, since
$$
\theta'_p \eta_{\Delta[p]}  = \pi\eta_{\Delta[p]} +
s\theta'_{p-1}g_p\eta_{\Delta[p]} - sf_p\pi\eta_{\Delta[p]} =
{\id} + sf_p - sf_p = \id,
$$
where we have used the commutativity proved above.
\end{proof}

\subsection{} Any contravariant monoidal functor $F^\ast :
(\Sets)^o\longrightarrow \Cochains{R}$ gives rise to a functor on
$\Sets$ with values in the category of $R$-differential graded
algebras by composing the K{\"u}nneth morphism of $F^\ast$ with
the morphism induced by the diagonal of the space
$$
F^\ast (X)\otimes F^\ast (X)\stackrel{\kappa}{\longrightarrow}
F^\ast (X\times X)\stackrel{\Delta^\ast}{\longrightarrow}
F^\ast(X).
$$
Now, from theorem \ref{unicidadderham} we deduce
\begin{corollary}
Let $A^\ast_\bullet$ be a simplicial $R$-dg algebra that satisfies
axioms (a)-(c). For any simplicial set $X$, the algebras
$S^\ast(X;R)$ and $A^\ast(X)$ are weakly equivalent in the
category of $R$-differential graded algebras.
\end{corollary}

\begin{example}
As application we obtain a proof, for any simplicial set $X$, of
the equivalence of the cochain algebra $S^\ast(X;\mathbb{Q})$ and
the polynomial De Rham algebra (compare with \cite{S} and
\cite{BG}). For this, take $L_{\bullet}^{\ast}$ the simplicial
cochain complex defined by the regular differential forms on the
$\mathbb Q$-cosimplicial scheme $H^\bullet$, whose $p$ component
is the hyperplane $H^p$ of $\mathbb{A}^{p+1}$ defined by the
equation $t_0 + \dots + t_p = 1$. It is a simplicial $\mathbb
Q$-dga that satisfies the homology and homotopy axioms (see
\cite{BG}), so if we define the algebra of polynomial De Rham
forms of $X$ as
$$
Su^\ast(X) = {\rm
Hom}_\Delta(X_\bullet,L^\ast_\bullet(\mathbb{Q})),
$$
the last corollary ensures that $S^\ast(X;\mathbb{Q})$ and
$Su^\ast(X)$ are weakly equivalent $\mathbb Q$-differential graded
algebras. The rational field $\mathbb Q$ may be replaced by any
field $\mk$ of characteristic zero.
\end{example}

\end{large}
\end{document}